\documentclass{article}

\usepackage[titletoc,toc,title]{appendix}

\usepackage{amssymb,amsfonts}
\usepackage{pifont}
\usepackage{amsthm}
\usepackage{amsmath}
\usepackage{makecell} 
\usepackage{multirow}
\usepackage{amsfonts}
\usepackage{amssymb} 
\usepackage{latexsym}
\usepackage{mathtools}
\DeclarePairedDelimiter{\ceil}{\lceil}{\rceil}
\usepackage{pifont}
\usepackage{graphics,graphicx,psfrag}
\usepackage{amscd}
\usepackage{epsfig}
\usepackage[all]{xy}
\usepackage{float}
\usepackage{enumitem}
\usepackage{cancel}
\usepackage[colorlinks=true]{hyperref}
\usepackage[hang,flushmargin]{footmisc} 
\usepackage[usenames,dvipsnames]{color}
\usepackage[usenames,dvipsnames,table]{xcolor}
\usepackage{tikz,siunitx}
\usetikzlibrary{matrix,arrows,positioning,calc}
\usetikzlibrary{decorations.pathmorphing}
\usepackage{MnSymbol}
\usepackage{setspace}
\usepackage{parskip}
\usepackage{rotating}

\let\OLDthebibliography\thebibliography
\renewcommand\thebibliography[1]{
  \OLDthebibliography{#1}
  \setlength{\parskip}{10pt}
  \setlength{\itemsep}{0pt plus 0.3ex}
}

\input xy
\xyoption{all} \xyoption{poly} \xyoption{knot}\xyoption{curve}
\input{diagxy}

\def\tn{\textnormal}

\newcommand{\cat}[1]{\mathcal{#1}}

\def\T{\mathbb T}

\renewcommand{\colon}{\nobreak\mskip2mu\mathpunct{}\nonscript
  \mkern-\thinmuskip{:}\mskip6muplus1mu\relax}

\def\Hom{\tn{Hom}}
\def\id{\tn{id}}

\def\Ob{\tn{Ob}}
\def\Mor{\tn{Mor}}

\def\ssec{{\color{blue!75!green}{\S}}}
\def\Defcolor{Definition}

\def\Figcolor{Figure}

\def\ie{{\em i.e.}, }
\def\eg{{\it e.g.}, }

\def\to{\rightarrow}

\def\inj{\hookrightarrow}

\def\too{\longrightarrow}

\def\ss{\subseteq}

\def\|{{\;|\;}}
\def\m1{{-1}}

\renewcommand{\ul}[1]{#1}

\newcommand{\qt}[1]{\tn{$\langle\!\langle#1\rangle\!\rangle$}}

\def\bhline{\Xhline{2\arrayrulewidth}}
\def\bbhline{\Xhline{2.5\arrayrulewidth}}

\newcommand{\boxtitle}[1]{\begin{center}#1\end{center}\vspace{-.1in}}

\def\ullimit{\ar@{}[rd]|(.3)*+{\lrcorner}}
\def\urlimit{\ar@{}[ld]|(.3)*+{\llcorner}}
\def\lllimit{\ar@{}[ru]|(.3)*+{\urcorner}}
\def\lrlimit{\ar@{}[lu]|(.3)*+{\lrcorner}}
\newcommand{\clabel}[1]{\ar@{}[rd]|(.5)*+{#1}}

\newcommand{\arr}[1]{\ar@<.5ex>[#1]\ar@<-.5ex>[#1]}
\newcommand{\arrr}[1]{\ar@<.7ex>[#1]\ar@<0ex>[#1]\ar@<-.7ex>[#1]}
\newcommand{\arrrr}[1]{\ar@<.9ex>[#1]\ar@<.3ex>[#1]\ar@<-.3ex>[#1]\ar@<-.9ex>[#1]}
\newcommand{\arrrrr}[1]{\ar@<1ex>[#1]\ar@<.5ex>[#1]\ar[#1]\ar@<-.5ex>[#1]\ar@<-1ex>[#1]}

\newcommand{\To}[1]{\xrightarrow{#1}}

\newcommand{\From}[1]{\xleftarrow{#1}}
\newcommand{\obox}[1]{\fbox{\mbox{{\color{white}p}\hspace{-1ex}#1{\color{white}h}\hspace{-1ex}}}}
\newcommand{\ovbox}[2]{\overset{#1}{\obox{#2}}}
\newcommand{\unbox}[2]{\underset{#1}{\obox{#2}}}
\newcommand{\fakebox}[1]{\tn{$\ceil{\mbox{#1}}$}}
\newcommand{\unit}{e}

\def\mfL{\mathcal{L}}

\def\N{\mathcal{N}}
\def\oN{\overline{N}}
\def\V{\mathcal{V}}
\def\oV{\overline{V}}
\def\L{\mathcal{L}}

\def\forget{U}
\def\oF{\overline{F}}
\def\omfL{\overline{L}}
\def\omfM{\overline{M}}

\def\Auth{\tn{Auth}}
\newcommand{\lozend}{~\hfill $\blacklozenge$}

\newcommand{\Adjoint}[4]{\xymatrix@1{#2 \ar@<.5ex>[r]^-{#1} & #3 \ar@<.5ex>[l]^-{#4}}}

\newtheoremstyle{monchito}{10pt}{6pt}{}{}{\color{blue!75!green}\bfseries}{\;\;}{ }{}

\theoremstyle{monchito}
\newtheorem{definition}{Definition}[subsection]

\theoremstyle{monchito}

\theoremstyle{monchito}
\newtheorem{proposition}[definition]{Proposition}

\theoremstyle{monchito}
\newtheorem{remark}[definition]{Remark}

\theoremstyle{monchito}

\theoremstyle{monchito}
\newtheorem{example}[definition]{Example}

\theoremstyle{monchito}

\theoremstyle{monchito}

\theoremstyle{monchito}

\theoremstyle{monchito}
\newtheorem{postulate}[definition]{Linguistic postulate}

\usepackage{newfloat}
\usepackage{caption}
\DeclareFloatingEnvironment[name=Figure]{myfigure}
\DeclareFloatingEnvironment[name=Figure]{myotherfigure}
\DeclareFloatingEnvironment[name=Figure]{anotherfigure}

\usepackage[margin=1.5in]{geometry}

\setlist[itemize]{parsep=\baselineskip}

\newcommand{\Set}{ \mathsf{Set} }                                     
\newcommand{\Olog}{\mathsf{Olog}}                                   
\newcommand{\InstOlog}{\mathsf{InstOlog}}                       
\newcommand{\Cat}{\mathsf{Cat}}                                       
\newcommand{\Eng}{\mathsf{Eng}}                                      
\newcommand{\iEng}{\mathsf{InstEng}}                               

\hypersetup{
    bookmarks=true,         
    unicode=false,          
    pdftoolbar=true,        
    pdfmenubar=true,        
    pdffitwindow=false,     
    pdfstartview={FitH},    
    pdftitle={My title},    
    pdfauthor={Author},     
    pdfsubject={Subject},   
    pdfcreator={Creator},   
    pdfproducer={Producer}, 
    pdfkeywords={keyword1} {key2} {key3}, 
    pdfnewwindow=true,      
    colorlinks=true,       
    linkcolor=blue!75!green,          
    citecolor=magenta,        
    filecolor=cyan,      
    urlcolor=violet           
}

\usepackage{titlesec}

\titleformat{\section}
{\color{blue!75!green}\normalfont\Large\bfseries}
{\color{blue!75!green}\thesection}{1em}{}

\titleformat{\subsection}
{\color{blue!75!green}\normalfont\large\bfseries}
{\color{blue!75!green}\thesubsection}{1em}{}

\titleformat{\subsubsection}
{\color{blue!75!green}\normalfont\large\bfseries}
{\color{blue!75!green}\thesubsubsection}{1em}{}

\tikzset{snake it/.style={-stealth,
decoration={snake, 
    amplitude = .4mm,
    segment length = 2mm,
    post length=0.9mm},decorate}}

    \AtBeginDocument{
  \DeclareSymbolFont{AMSb}{U}{msb}{m}{n}
  \DeclareSymbolFontAlphabet{\mathbb}{AMSb}}

\usepackage{microtype}

\begin{document}

\title{{\bf Toward formalizing ologs:} \\ {\bf {\normalsize linguistic structures, instantiations, and mappings}}}

\author{{\sc Marco A. P\'erez} \\ 
{\footnotesize Department of Mathematics} \\ 
{\footnotesize Massachusetts Institute of Technology} \\ 
{\footnotesize \texttt{maperez$@$mit.edu}} 
\and 
{\sc David I. Spivak} \\ 
{\footnotesize Department of Mathematics} \\ 
{\footnotesize Massachusetts Institute of Technology} \\ 
{\footnotesize \texttt{dspivak$@$math.mit.edu}} 
}

\date{\today}

\maketitle

\begin{abstract} \setstretch{0.75} \noindent \footnotesize We define the notion of linguistic structure on a small category, in order to provide a more formal description of ontology logs, also known as ologs, introduced in \cite{SpivakOlogs} by R. E. Kent and the second author. In particular, we construct a bicategory $\Eng$, of English noun phrases and verb phrases, endorsed as functional by varying sets of authors. An olog is then defined as a lax functor to $\Eng$. We then present a new notion of linguistic functor, which extends the notion of meaningful functors defined in \cite{SpivakOlogsPre}. Finally, we discuss the relationship between ologs and databases in this context. 
\end{abstract}

\setcounter{tocdepth}{2}
\tableofcontents


\section{Introduction}\label{sec:intro}

The theory of ontology logs (ologs for short) was introduced by Robert Kent and the second author in their paper \cite{SpivakOlogs}, as a framework for knowledge representation. Ologs are basically mathematical categories that have been wrapped in natural-language English. They have been applied in several branches of science and engineering, \cite{Bergstra,Cranford,Giesa,Gomez,Kent,SpivakMaterial}, as a tool for various kinds of formal modeling.

Typically, a person who wishes to record and document some of her knowledge or ideas will do so in prose, \eg a scientist publishes ideas in the form of research papers. Ologs offer the ability to express complex ideas using a special type of diagrams. Namely, the objects of study and the relationships between them can be represented as the objects and arrows in a category. The difference between an olog and a category is that an olog has additional structure: each object is labeled with a noun phrase and each arrow is labeled with a verb phrase, so that reading source-arrow-target yields an English sentence. These must satisfy certain rules, giving a set-theoretic semantics to the olog, which in turn allows it to serve a dual role as a database schema. There is a formula for composing sentences end-to-end into a new sentence, when following a path of arrows through an olog, and a pair of equivalent paths (also known as a commutative diagram) in the category is understood as a \textit{declared fact} equating the two English sentences.    

The idea of regarding English sentences as mappings is not new. One interesting approach on this matter is given by the notion of conceptual metaphor \cite{WikiConceptual} which links one idea (a \textit{source domain}) to another (a \textit{target domain}) to better understanding something. More formal approaches (compared to the previous one) to language, knowledge and information modeling have been developed by other authors, such as D. Kartsaklis, M. Sadrzadeh, S. Pulman and B. Coecke. In \cite{Coecke}, for instance, they make use of notions from category theory, such as \textit{compact closed categories} and \textit{strongly monoidal functors}, to studying meaning in natural language. 

The primary goal of this paper is to present a reformulation of the previous concept of olog, which leads to a categorical and linguistic description of mappings between ologs. In our reformulation, any type, aspect, or fact (in the sense of \cite{SpivakOlogs}) in an olog must be endorsed by a set of people, who understand its meaning, \ie the way it indicates a set or a function. This cultural understanding of the English language is captured by a certain bicategory $\Eng$. The above idea that ologs are categories wrapped in English will be formalized by saying that ologs are categories mapping to $\Eng$. Our notion of mapping between ologs falls out of that structure. The result is that a map between ologs is not just a functor between their underlying categories; it is a functor that respects the linguistic description on each node and arrow of the ologs. For example, there is an obvious functor $F$ (namely the identity) between the underlying categories of the following ologs:

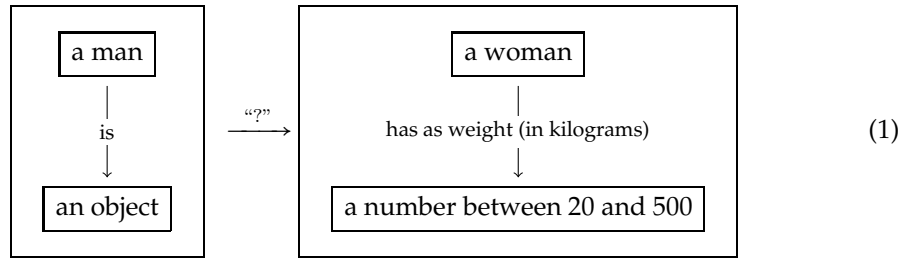
\begin{equation}\label{eqn:no_chance} 
\parbox{1.1in}{\fbox{ \begin{tikzpicture}[description/.style={fill=white,inner sep=2pt}] 
\matrix (m) [ampersand replacement=\&, matrix of math nodes, row sep=3.5em, column sep=0em] 
{ \obox{a man} \\ \obox{an object} \\ }; 
\path[->] 
(m-1-1) edge node[description] {$\footnotesize\mbox{is}$} (m-2-1); 
\end{tikzpicture} } 
}
\xrightarrow{\;\;``?"\;\;}
\parbox{2.3in}{\fbox{ \begin{tikzpicture}[description/.style={fill=white,inner sep=2pt}] 
\matrix (m) [ampersand replacement=\&, matrix of math nodes, row sep=3.5em, column sep=0em] 
{ \obox{a woman} \\ \obox{a number between 20 and 500} \\ }; 
\path[->] 
(m-1-1) edge node[description] {$\footnotesize \mbox{has as weight (in kilograms)}$} (m-2-1) ; 
\end{tikzpicture} }
}
\end{equation} 

But would this functor $F$ \emph{mean} anything? We will rephrase this question in \ssec~\ref{sec:com} as follows: ``does there exist an author who is willing to endorse a \textit{linguistic structure} on $F$?" In this paper we explain the difference between a (linguistic) map of ologs and a mere functor between their underlying categories. We will address this particular case (\ref{eqn:no_chance}) in Example~\ref{ex:no_chance}. 

To some extent, this issue is considered in \cite{SpivakOlogsPre}, where the authors introduce the concept of {\it meaningful functor}. One limitation of their approach is that it depends on a strong assumption, namely that every olog $\cat{C}$ is equipped with a functor $I \colon \cat{C} \too \Set$, and that this functor somehow controls the meaning of the olog. Such a set-valued functor is called an \textit{instantiation}, a term coming from database theory \cite{SpivakData}. The idea is that $I$ represents a kind of database of examples, or instances, for the various types, aspects, and facts in the olog. For example, if an author is writing an olog $\cat{C}$ describing a familiar real-world situation, then for a type $c$ (\eg $c=$ cat) the set $I(c)$ represents all the examples of $c$ (all cats) known by the author. Somewhat strangely, however, there is no requirement in \cite{SpivakOlogsPre} that the database functor $I$ should in any way correspond to the linguistic structure on the olog. A similar issue exists for functors between ologs in \cite{SpivakOlogsPre}. 

In this paper, we remedy these issues. First, we allow ologs to exist without being instantiated; that is we disentangle ologs and their instantiations. This way, the set of documented examples can evolve over time, without changing the olog to which they refer. On the other hand, we add a constraint to instantiations: for a set-valued functor to count as an instantiation of an olog $\cat{C}$, it must conform to the linguistic structure, the labelings, on $\cat{C}$. The same goes for mappings between ologs: in order for a functor to count as a mapping between ologs, it must conform to the linguistic structures involved. We also take more care to explain the relationship between an olog and its authors. We introduce the concept of \textit{endorsement}: an author can endorse that a certain concept or relationship between concepts makes sense, that a certain fact is true, etc. 

Here is an expert-level view of this paper. There is a bicategory $\Eng$ that denotes the English language as it divides into noun phrases, which indicate sets, and verb phrases, which indicate functions, where all of this ``indicating" is decided solely by speakers. Given that such a bicategory $\Eng$ exists, an olog is just a small category $\cat{C}$ and a lax functor $L \colon \cat{C} \too \Eng$, called a linguistic structure on $\cat{C}$. Allowing the base category to vary, we get a fibration $\Olog \too \Cat$, where $\Olog$ denotes the category of ologs and $\Cat$ is the category of small categories. Instantiated ologs are defined similarly: there is a bicategory $\iEng\ss\Eng\times\Set$ in which the noun phrase associated to each object is exemplified by a set, and each verb phrase associated to a morphism is exemplified by a function. An instantiated olog is a lax functor $\cat{C} \too \iEng$, and we again have a fibration $\InstOlog \too \Cat$. All of this will be explained in the main sections of the paper.

This paper is organized as follows: In \S~\ref{sec:EngBicat} we introduce the bicategory $\Eng$ of English language. We present the notion of \textit{author endorsement} for noun phrases, verb phrases and equivalence between sentences, along with a list of linguistic guidelines needed to design any olog within $\Eng$. \ssec~\ref{sec:def} introduces the category $\Olog$ of ologs, as well as the category $\InstOlog$ of instantiated ologs, by defining mappings between ologs and between instantiated ologs. To do so, we introduce the notion of \emph{linguistic functors} and \emph{instantiated functors}, as mentioned above. The latter of these is an adaptation of the ``meaningful functor" notion defined in \cite{SpivakOlogsPre}. 

Most of the assertions on the bicategory $\Eng$ have a linguistic element and so are not purely mathematical. However, with the help of several linguistic postulates (mainly found in \ssec~\ref{sec:endoresement_postulates}), fairly formal proofs are possible. Once $\Eng$ is given, the rest is straightforward category theory, thus giving to the theory of ologs a more solid theoretical basis. As the paper is mainly written for a non-mathematical audience, all proofs are given as prose arguments, merged with the text.


\subsection*{Background and Notation} 

We will assume the reader is familiar with some basic concepts from category theory, such as {\it opposite categories}, {\it isomorphisms}, {\it functors}, and  {\it natural transformations}. Readers without category-theoretic background may still benefit from reading the less categorical definitions and results, skimming the category theory, and trying to digest the examples. The books \cite{Adamek, SpivakBook} are good sources for category theory, with many illustrations. No previous knowledge on bicategories, lax functors and lax transformations is needed; everything we say about these topics will be spelled out in concrete terms, but interested readers can check the book \cite[Chapter 1]{Leinster} as an excellent source for a brief review on this matter. 

The word ``category" will always mean a small category unless otherwise stated, \ie the collections of objects and morphisms (also called arrows) are sets. Given a category $\cat{C}$, we will denote the sets of objects and morphisms of $\cat{C}$ by $\Ob(\cat{C})$ and $\Mor(\cat{C})$, respectively. For every two objects $c, c' \in \Ob(\cat{C})$, we will denote the set of morphisms from $c$ to $c'$ by $\Hom_{\cat{C}}(c,c')$. The composition of two morphism $f_1 \colon c \to c'$ and $f_2 \colon c' \to c''$ will be denoted using the semicolon symbol ``;", \eg by $f_1; f_2 \colon c \to c''$. Our reason to use this notation is to make the composition on morphisms in $\cat{C}$ parallel to the concatenation of sentences (see Definition~\ref{def:Concatenation}) in any olog set on $\cat{C}$. We denote the (large) category of sets, functions, and function composition by $\Set$.


\section{English language as a bicategory}\label{sec:EngBicat}

One of the goals of this paper is to formalize the definition of olog. As we mentioned in the introduction, the idea is to define an olog as a category $\cat{C}$ together with a lax functor from $\cat{C}$ to a particular bicategory $\Eng$, which encodes how English speakers understand their language. Details about bicategories, lax functors, and lax transformations will be given in the context of ologs, through the constructions of $\Eng$, linguistic structures, and linguistic functors, respectively. 

We begin this section presenting the notions of noun phrase, verb phrase, and equivalence of sentences in \S~\ref{sec:sentences}. These, along with the notion of author endorsement are going to be principle components defining the bicategory $\Eng$, presented in \S~\ref{sec:Eng}. As ologs are straddling the world of categories and linguistics, we must explain what it means for the labels on objects and morphisms to correspond to valid sets and functions. This concept of endorsement should satisfy certain guidelines, which we will present as (linguistic) postulates. After assuming these postulates, we will be able in \S~\ref{sec:Eng} to regard the English language as a bicategory $\Eng$, whose collections of $0$-cells, $1$-cells, and $2$-cells formalize the concepts of types, aspects, and facts (respectively) that define an olog in the spirit of \cite{SpivakOlogs}. We formalize ologs as lax functors from small categories to $\Eng$ (see Definition~\ref{def:Linguistic}). 

In \S~\ref{sec:Eng_instances} we will equip each type in $\Eng$ with a set of instances, \ie an author-specified set of examples corresponding to the linguistic label. Doing so, we obtain another bicategory, denoted $\iEng$, and following the same plan as above, we define instantiated ologs later in Definition~\ref{def:Instantiation}.


\subsection{Nouns, verbs, and functional sentences}\label{sec:sentences}

When we speak of a noun phrase, we refer to what is more formally called \emph{a singular indefinite noun phrase}, such as ``a person". If $\N$ is a symbol representing this noun phrase, we write \qt{\N}=``a person"; we call $\qt{\N}$ the {\bf reading of $\N$} (and say that $\N$ {\bf is read} \qt{\N}). We use similar notation for verb phrases. For us, a verb phrase $\V$ is that which connects one (singular indefinite) noun phrase, often called \emph{the subject}, to another, often called \emph{the object}. If \qt{\V}=``is" and \qt{\N'}=``a mammal", then we denote the concatenated string, ``a person is a mammal" by \qt{\N}\qt{\V}\qt{\N'} or simply by \qt{\N \V \N'}.

\begin{definition}\label{def:NounSet}
A \textbf{noun phrase} is formed by the indefinite article \emph{a/an} followed by an English noun which refers to a set-like concept, \ie a concept for which an author can record a set of examples. ~%
\lozend
\end{definition}

\begin{example}
The noun phrase \fakebox{a shirt}~%
\footnote{As in \cite[\S 2.1]{SpivakOlogs}, we write $\ceil{\mbox{a shirt}}$ instead of $\fbox{\mbox{a shirt}}$, because some typographical problems emerge when writing a text-box in a line of text. Note that the text-box seems out of place in this paragraph, and so many in-line text-boxes are troublesome for the aesthetic of this paper.} 
refers to the set of shirts, \ie things that can be called ``a shirt". Similarly, the noun phrase \fakebox{a prime number less than 10} refers to the set $\{2,3,5,7\}$.~ 
\lozend
\end{example}

In order to clarify the semantics, we impose a certain rule on our verb phrases: they must refer to functions, in the sense of set theory (see \cite{SpivakOlogs}).

\begin{definition}\label{def:FunctionRel} 
Given two noun phrases $\N_1$ and $\N_2$ and a verb phrase $\V$, one says that $\V$ \textbf{functionally connects} $\N_1$ to $\N_2$ if the concatenation $\qt{\N_1 \V \N_2}$ is an English sentence that refers to a mathematical function. In this case one says that $\Sigma = (\N_1, \V, \N_2)$ is a \textbf{functional sentence}.~ 
\lozend
\end{definition}

For example, given two noun phrases $\ceil{\mbox{a person}}$ and $\ceil{\mbox{a woman}}$, the arrow  
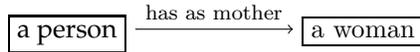
\[ \begin{tikzpicture} 
\matrix (m) [matrix of math nodes, row sep=0em, column sep=6em, text height=1.5ex, text depth=0.25ex]
{ \fbox{\mbox{a person}} & \fbox{\mbox{{\rm a woman}}} \\ };
\path[->]
(m-1-1) edge node[above] {$\mbox{{\rm \footnotesize has as mother}}$} (m-1-2);
\end{tikzpicture} \]
can be read as a sentence, ``a person has as mother a woman". This sentence expresses that for anything that could be called ``a person" there is something it ``has as mother" that can itself be called ``a woman". This seems true, but more important for us, it seems to represent a function: each person has only one mother. A reader who understands these concepts probably has an example of a person (namely himself) and the woman that corresponds to the mother he has. Thus the arrow text is a verb phrase that functionally connects $\ceil{\mbox{a person}}$ and $\ceil{\mbox{a mother}}$
in the sense of Definition~\ref{def:FunctionRel}. 

We can contrast the above situation with the following:
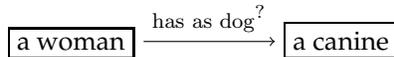
\begin{equation*}
\begin{tikzpicture} 
\matrix (m) [matrix of math nodes, row sep=0em, column sep=5em, text height=1.5ex, text depth=0.25ex]
{ \fbox{\mbox{a woman}} & \fbox{\mbox{a canine}} \\ };
\path[->]
(m-1-1) edge node[above] {$\mbox{{\rm \footnotesize has as dog}}^?$} (m-1-2);
\end{tikzpicture}
\end{equation*}
Because not every woman has a dog, and some women have two dogs, the arrow text is a verb phrase that \emph{does not} correspond to a function; \ie it does not functionally connect its source and target noun phrases. 

We next consider the notion of equivalent verb phrases, by which we mean verb phrases that indicate the same functional relationship. For example, if I say ``An integer has as successor an integer", and you say ``An integer yields, by adding $1$, an integer", we are saying the same thing using different verbs.

\begin{definition}\label{def:DeclaredEq} 
Two sentences $\Sigma = (\N_1, \V, \N_2)$ and $\Sigma' = (\N_1, \V', \N_2)$ are said to be \textbf{equivalent}, denoted $\Sigma \simeq \Sigma'$, if they refer to the same functional relationship. If $\Sigma \simeq \Sigma'$, then their verb phrases are equivalent as well, and we denote $\V \simeq \V'$.  

Graphically, we will draw a checkmark $\checkmark$ between the two sentences $\Sigma$ and $\Sigma'$ if they---and hence their verb phrases $\V, \V'$---are equivalent, as shown in the following figure. 
\begin{figure}[H]
\centering
\belowcaptionskip = -15pt
\begin{tikzpicture}[description/.style={fill=white,inner sep=2pt}]
\matrix (m) [matrix of math nodes, row sep=1.15em, column sep=2em, text height=1.75ex, text depth=1.25ex]
{ \N_1 &\checkmark & \N_2 \\ };
\path[->]
(m-1-1) edge [bend left=20] node[above] {\footnotesize$\V$} (m-1-3)
(m-1-1) edge [bend right=20] node[below] {\footnotesize$\V'$} (m-1-3);
\end{tikzpicture} 
\caption{Equivalence between $\V$ and $\V'$.}\label{fig:equivalence_V_Vprime}
\end{figure}
As in \cite[\S 2.3.3.4]{SpivakBook}, an equivalence $\Sigma \simeq \Sigma'$ can be given the following English-language interpretation (where $x, y_1,$ and $y_2$ are just symbols, to be copied verbatim):
\begin{figure}[H]
\centering
\belowcaptionskip = -15pt
\begin{quote}
{\color{red!40!black}{For any $\qt{\N_1}$ $x$, we know that $x$ $\qt{\V \N_2}$, that we call $y_1$, and we know that $x$ $\qt{\V' \N_2}$, that we call $y_2$; and the fact is, $y_1$ and $y_2$ are the same for any $x$.}}
\end{quote}
\caption{How equivalences are interpreted as English.}
\label{fig:FormulaFact}
\end{figure}~%
\lozend
\end{definition}

\begin{example} 
We explain how to use the syntax formula in Figure~\ref{fig:FormulaFact} to show that the following two sentences are equivalent:
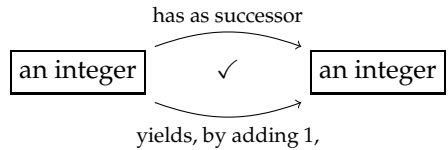
\begin{center}
\begin{tikzpicture}[description/.style={fill=white,inner sep=2pt}]
\matrix (m) [matrix of math nodes, row sep=1em, column sep=2em, text height=1.75ex, text depth=1.25ex]
{\fbox{\tn{an integer}} & \checkmark& \fbox{\tn{an integer}} \\ };
\path[->]
(m-1-1) edge [bend left=20] node[above] {\footnotesize has as successor} (m-1-3)
(m-1-1) edge [bend right=20] node[below] {\footnotesize yields, by adding 1,} (m-1-3);
\end{tikzpicture} 
\end{center}
The equivalence $\Sigma \simeq \Sigma'$ is read as follows:
\begin{quote}
For any integer $x$, we know that $x$ has as successor an integer, that we call $y_1$, and we know that $x$ yields by adding $1$ an integer, that we call $y_2$; and the fact is, $y_1$ and $y_2$ are the same for any $x$.
\end{quote}
We can agree that the previous statement is valid, since $y_1 = x + 1 = y_2$, by the definition of successor, and thus $\Sigma \simeq \Sigma'$. ~%
\lozend
\end{example}

\begin{definition}\label{def:Concatenation} 
Let $\Sigma_1 = (\N_1, \V_1, \N_2)$ and $\Sigma_2 = (\N_2', \V_2, \N_3)$ be two sentences. If $\N_2 = \N_2'$, we say that $\Sigma_1$ and $\Sigma_2$ are \textbf{concatenatable}. Let $\V_1;\V_2$ represent the verb phrase
$$
\qt{\V_1;\V_2}\coloneqq \qt{\V_1 \N_2} \mbox{``which"} \qt{\V_2}.
$$
We define the \textbf{concatenation} of $\Sigma_1$ and $\Sigma_2$, denoted $\Sigma_1;\Sigma_2$, to be the sentence\\\vspace{-.2in}

\centering$\Sigma_1;\Sigma_2 := (\N_1, \V_1;\V_2, \N_3).$\\
\vspace{-.4cm}\hfill$\blacklozenge$
\end{definition}

\begin{example}\label{ex:amino_composite}

Consider the diagram below:
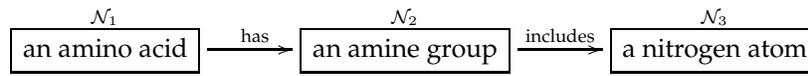
\[
\xymatrix@C=32pt{
\ovbox{\N_1}{an amino acid}\ar[r]^{\tn{has}}&\ovbox{\N_2}{an amine group}\ar[r]^{\tn{includes}}&\ovbox{\N_3}{a nitrogen atom}
}
\]
If $\qt{\V_1} = \mbox{``has"}$ and $\qt{\V_2} = \mbox{``includes"}$, then the sentences $\Sigma_1=(\N_1,\V_1,\N_2)$ and $\Sigma_2=(\N_2,\V_2,\N_3)$ are concatenatable. Their concatenation is read
\begin{equation}\label{eqn:concatenation_example}
\qt{\Sigma_1;\Sigma_2} = \tn{``an amino acid has an amine group, which includes a nitrogen atom".} 
\end{equation}
~%
\lozend
\end{example}

\begin{example}\label{ex:ExampleFacts} 
In this example, we combine the concepts of equivalence and concatenation of sentences. The following diagram represents an equivalence between sentences, as indicated by the checkmark symbol $\checkmark$: 
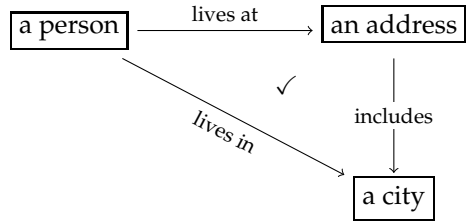
\[ \begin{tikzpicture}[description/.style={fill=white,inner sep=2pt}]
\matrix (m) [matrix of math nodes, row sep=1.15em, column sep=3em, text height=1.75ex, text depth=1.25ex]
{ \fbox{\mbox{a person}} & & \fbox{\mbox{an address}} \\ & \mbox{} \\ & & \fbox{\mbox{a city}} \\ };
\path[->]
(m-1-1) edge node[above] {$\footnotesize \mbox{lives at}$} (m-1-3) edge node[below,sloped] {$\footnotesize \mbox{lives in}$} (m-3-3)
(m-1-3) edge node[description] {$\footnotesize \mbox{includes}$} (m-3-3)
(m-1-3)-- node[description] {$\checkmark$} (m-2-2);
\end{tikzpicture} \]
First, note that the concatenation of the top and right arrows yields the sentence ``a person lives at an address, which includes a city", and the diagonal arrow is read ``a person lives in a city". According to \Figcolor~\ref{fig:FormulaFact}, the equivalence between these two sentences is read as the assertion
\begin{quote}
For any person $x$, we know that $x$ lives at an address which includes a city $y_1$, and we know that $x$ lives in a city $y_2$; and the fact is, $y_1$ and $y_2$ are the same for any $x$.
\end{quote}
Contrast that fact to the following diagram:

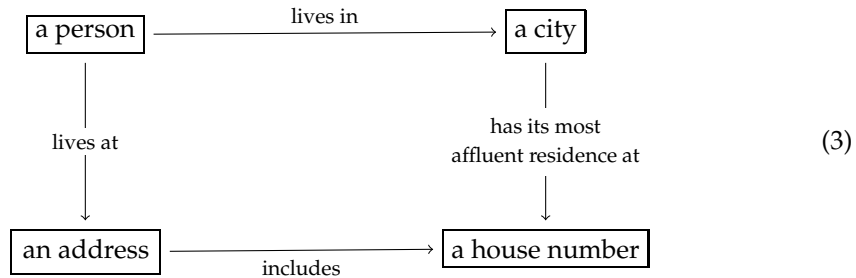
\begin{equation}\label{eqn:non-equiv} 
\parbox{3.5in}{
\begin{tikzpicture}[description/.style={fill=white,inner sep=2pt}]
\matrix (m) [ampersand replacement=\&, matrix of math nodes, row sep=6em, column sep=10em]
{ \fbox{\mbox{a person}} \& \fbox{\mbox{a city}} \\ \obox{an address} \& \fbox{\mbox{a house number}} \\ };
\path[->]
(m-1-1) edge node[above] {$\footnotesize \mbox{lives in}$} (m-1-2) edge node[description] {$\footnotesize \mbox{lives at}$} (m-2-1)
(m-1-2) edge node[description] {$\footnotesize \begin{array}{c} \mbox{has its most} \\ \mbox{affluent residence at} \end{array}$} (m-2-2)
(m-2-1) edge node[below] {$\footnotesize \mbox{includes}$} (m-2-2);
\end{tikzpicture}
}
\end{equation}
which \emph{does not} represent a valid equivalence (note the absence of the checkmark symbol $\checkmark$). In this case, the assertion indicating the equivalence between the two involved sentences would be read
\begin{quote}
For any person $x$, we know that $x$ lives in a city, which has its most affluent residence at a house number $y_1$, and we know that $x$ lives at an address, which includes a house number $y_2$; and the fact is, $y_1$ and $y_2$ are the same for any $x$. 
\end{quote} 
This assertion is quite dubious, because there seem to be many city-dwellers whose house numbers are different from that of the most affluent resident in their city. Hence, we as authors cannot endorse \eqref{eqn:non-equiv} as an equivalence. ~%
\lozend
\end{example}

\begin{definition}\label{def:linguistic_expression}
We say that a symbol $\L$ is a \textbf{linguistic expression} if it represents a noun phrase, a verb phrase, a sentence, or an equivalence between sentences. ~%
\lozend
\end{definition}


\subsection{Some postulates about endorsement}\label{sec:endoresement_postulates}

In this section, we introduce the notion of author endorsement. This is the cultural aspect of the English (or any) language: noun phrases and verb phrases do not inherently represent sets and functions. The relationship between words and their meaning is, to some extent, author dependent. 

After defining endorsement, we provide some postulates about it. These postulates are what allows us to represent the English language as a bicategory $\Eng$ in \S~\ref{sec:Eng}, which is the main ingredient in our formal definition of ologs. 

The following definition cannot, by its nature, be made mathematically rigorous. Indeed, it asks the reader to grant that there are things called people, and that people can evaluate whether verbs refer to functions, etc.

\begin{definition}\label{def:EndorsementNoun} 
Let $s$ be a person who understands Definitions~\ref{def:NounSet},~\ref{def:FunctionRel},~and~\ref{def:DeclaredEq}.
\begin{itemize}[topsep=-1ex,itemsep=-1.5ex]
\item[\textbf{(a)}] If $\N$ is a noun phrase, we will say that $s$ \textbf{endorses} $\N$, denoted $s \models \N$, if $\qt{\N}$ refers to a distinction made and recognizable by $s$, \ie if $s$ agrees that $\N$ is a set-like noun phrase, in the sense of Definition~\ref{def:NounSet}, and can imagine example elements of it.

\item[\textbf{(b)}] Suppose that $\Sigma = (\N_1, \V, \N_2)$ is a sentence and that $s\models\N_1$ and $s\models\N_2$. We say that $s$ \textbf{endorses} $\V$, denoted $s \models \V$, if $s$ is willing to declare that $\V$ is functional, in the sense of Definition~\ref{def:FunctionRel}, and understands how examples of $\N_1$ correspond via $\V$ to examples of $\N_2$. In this case we also write $s\models \Sigma$.

\item[\textbf{(c)}] If $\Sigma= (\N_1, \V, \N_2)$ and $\Sigma'= (\N_1, \V', \N_2)$ are two sentences such that $s \models \Sigma$ and $s \models \Sigma'$, then we say that $s$ \textbf{endorses} $\Sigma \simeq \Sigma'$, denoted $s \models [\Sigma \simeq \Sigma']$, if $s$ is willing to declare that $\Sigma \simeq \Sigma'$ is an equivalence in the sense of Definition~\ref{def:DeclaredEq} and that the correspondences $\V$ and $\V'$ act identically on examples of $\N_1$, as in \textbf{(b)}.  
\end{itemize}
By an \textbf{author set} $S$ we mean a set of people in which each member is able to decide whether or not a linguistic expression (see Definition~\ref{def:linguistic_expression}) is valid in the sense of \textbf{(a), (b), (c)} above. If $S$ is an author set and $\L$ is a linguistic expression, we say that $S$ \textbf{endorses} $\L$, denoted $S \models \L$, if $s \models \L$ for every $s \in S$. ~%
\lozend
\end{definition}

\begin{example}
A reader may not be able to endorse the sentences $\Sigma_1$ and $\Sigma_2$ of Example~\ref{ex:amino_composite} as functional, whereas others would endorse them. That is, some may not know whether every amino acid always contains exactly one amine group or whether every amine group contains exactly one nitrogen atom. But given that an author $s$ endorses them (and we should mention that $\Sigma_1$ and $\Sigma_2$ are broadly endorsed), Postulate~\ref{rules:EndorsementConcatenation} below says that $s$ must also endorse their concatenation $\Sigma_1;\Sigma_2$.\lozend
\end{example}

Recall that a postulate is an idea suggested or assumed as true as the basis for reasoning, discussion, or belief. We make a few linguistic postulates throughout this section and \S~\ref{sec:Eng_instances}, starting with the next one.

\begin{postulate}\label{rules:EndorsementConcatenation} 
Let $\Sigma_1$ and $\Sigma_2$ be two concatenatable sentences as in Definition~\ref{def:Concatenation}. If $s$ endorses both sentences $s \models \Sigma_1$ and $s \models \Sigma_2$, then we will assume that $s$ endorses their concatenation, $s \models \Sigma_1; \Sigma_2$. ~%
\lozend  
\end{postulate}

The following postulate says that every olog author must endorse sentences like ``a bottle is of course a bottle" as functional. We will later postulate (in~\ref{rules:EndorsementFact}) that such sentences correspond to identity functions on sets, in what we call \emph{instantiated English}.

\begin{postulate}\label{rules:I}
We will assume that there is a unique verb phrase $e$, read $\qt{e} = \mbox{``is of course"}$, such that if $\N$ is a noun phrase and $s$ is an author with $s\models\N$, then $s \models (\N, e, \N)$. We call $e$ the \textbf{unit verb phrase}. ~%
\lozend
\end{postulate}

The following postulate says that the unit verb phrase introduced above is unital with respect to the string concatenation defined in Definition~\ref{def:Concatenation}.

\begin{postulate}\label{rules:I_unital} 
Let $\unit$ be the unit verb phrase as in Postulate~\ref{rules:I}. For every sentence $\Sigma = (\N_1, \V, \N_2)$ and any author $s$, if $s \models \Sigma$ then~\\\vspace{-.2in}

\centering $s \models [\Sigma; (\N_2, \unit, \N_2) \simeq \Sigma]$\quad and\quad $s \models [(\N_1, \unit, \N_1); \Sigma \simeq \Sigma]$.\\
\vspace{-.4cm}\hfill$\blacklozenge$
\end{postulate}

\begin{remark}\label{rem:concatenation_associative} 
Associativity of concatenation is always endorsed by any author endorsing the sentences involved, but this does not require a postulate. It follows from the associativity of string concatenation. ~%
\lozend
\end{remark}


\subsection{The bicategory $\Eng$ of English expressions}\label{sec:Eng}

In this section we explain how to combine the notions of noun phrase, verb phrase, equivalence, and endorsement to construct a bicategory based on the English language.

Recall from Definition~\ref{def:EndorsementNoun} the notion of endorsement (denoted by the symbol $\models$) for noun phrases, verb phrases, and equivalence between sentences. Recall also that by an \textbf{author}, we mean a person who understands Definitions~\ref{def:NounSet},~\ref{def:FunctionRel},~and~\ref{def:DeclaredEq}; an \textbf{author set} is a set of authors.

\begin{definition}\label{def:type_aspect_fact} \
\begin{itemize}[topsep=-1ex,itemsep=-1.5ex]
\item[\textbf{(a)}] A \textbf{type} $N=(\N,T)$ consists of a noun phrase $\N$ (see Definition~\ref{def:NounSet}) and an author set $T$ such that $T \models \N$. The author set may be denoted $\Auth(N)\coloneqq T$.

\item[\textbf{(b)}] Given two types $N_1$ and $N_2$, an \textbf{aspect} from $N_1$ to $N_2$, denoted $V \colon N_1 \to N_2$, consists of a pair $V=(\V,A)$, where $\V$ is a verb phrase functionally connecting $\N_1$ to $\N_2$ in the sense of Definition~\ref{def:FunctionRel}, where $A$ is an author set such that $A \subseteq \Auth(N_1) \cap \Auth(N_2)$, and where $A \models \V$ . The author set may be denoted $\Auth(V)\coloneqq A$. 

\item[\textbf{(c)}] Given two types $N_1$ and $N_2$, and two aspects $V_1, V_2 \colon N_1 \to N_2$, a \textbf{fact} from $V_1$ to $V_2$, denoted $V_1 \Rightarrow V_2$, is given by an author set $F \subseteq \Auth(V_1) \cap\Auth(V_2)$ and an endorsement $F \models [\V_1 \simeq \V_2]$, in the sense of Definition~\ref{def:DeclaredEq}. We allow an abuse of notation and write $F \colon V_1 \Rightarrow V_2$.\footnote{Notice that we use ``$\checkmark$" only to denote equivalence between sentences, and the double arrow ``$\Rightarrow$" to denote $2$-cells in $\Eng$, which are formed by an equivalence between sentences \emph{and} an endorsing author set.} ~%
\lozend
\end{itemize}
\end{definition}

We will soon show (see Proposition~\ref{prop:Eng}) that \textbf{(a)}, \textbf{(b)}, and \textbf{(c)} above define the collection of $0$-cells, $1$-cells and $2$-cells of a bicategory, which we denote $\Eng$. Roughly speaking, we need to define two operations on $1$-cells and $2$-cells, usually called \textit{horizontal} and \textit{vertical compositions}, satisfying certain properties and \textit{coherence conditions}. 

Most of the work is dealing with author sets, a consideration with which many readers need not be overly concerned. The main idea, which one should really know in order to proceed, is that ordinary (``horizontal") composition in $\Eng$ involves the concatenation of sentences, from Definition~\ref{def:Concatenation}. On a first reading, one may skip to \S~\ref{sec:Eng_instances}. On a more careful reading, on the other hand, a good way to understand this point in a more formal way, is to study first the collection of all noun phrases and verb phrases endorsed by a single author $s$, which will form a category $\Eng_s$. In that particular context, author endorsement becomes straightforward, allowing us to treat the composition of morphisms in $\Eng_s$ as a mere concatenation of sentences. We provide the details about this in the next section. Once we have defined $\Eng_s$, we will focus on proving Proposition~\ref{prop:Eng}.


\subsubsection*{The category of English expressions endorsed by a single author}\label{sec:Engs}

\begin{proposition}\label{prop:Engs} 
Let $s$ be an author. There is a category $\Eng_s$ with the following collection of objects and morphisms:
\begin{itemize}[topsep=-1ex,itemsep=-1.5ex]
\item[\textbf{(a)}] $\Ob(\Eng_s)$ is the set of types in $\Eng$ which have the form $N = (\N, s)$, \ie $s \models \N$.

\item[\textbf{(b)}] $\Mor(\Eng_s)$ is the set of aspects in $\Eng$ which have the form $V = (\V, s)$, \ie $s \models \V$.  ~%
\lozend
\end{itemize}
\end{proposition}

Proposition~\ref{prop:Engs} is basically a consequence of some of the linguistic postulates in \S~\ref{sec:endoresement_postulates}. First, given two composable morphisms $V_1 \colon N_1 \to N_2$ and $V_2 \colon N_2 \to N_3$, the composition of $V_1$ and $V_2$ in $\Eng_s$ is given by the concatenated verb phrase $\V_1; \V_2$, which is endorsed by $s$ by Postulate~\ref{rules:EndorsementConcatenation}, and which is associative by Remark~\ref{rem:concatenation_associative}. Finally, for every object $N$ in $\Eng_s$, one has that $s$ endorses the sentence $(\N, e, \N)$ by Postulate~\ref{rules:I}, because $s \models \N$. The morphisms of the form $(e,s) \colon N \to N$, which we denote by ${\rm id}_N$, are the identity morphisms in $\Eng_s$ by Postulate~\ref{rules:I_unital}.

One interesting property about $\Eng_s$ is that the equivalence between sentences is a \textit{congruence relation} $\simeq_s$ on $\Eng_s$, meaning it satisfies the conditions of Proposition~\ref{prop:congruence}. For every $N, N' \in \Ob(\Eng_s)$ and every $V, V' \in \Hom_{\Eng_s}(N,N')$, we write
 $$V \simeq^{N,N'}_s V'$$ 
 if $s \models [V \simeq V']$.
The following result provides the necessary theoretical framework to define the \textit{horizontal composition} of aspects and the \textit{vertical composition} of facts in the bicategory structure of $\Eng$.

\begin{proposition}\label{prop:congruence} For every author $s$, the relation $\simeq_s$ satisfies the following properties:
\begin{itemize}[topsep=-1ex,itemsep=-1.5ex]
\item For every $N, N' \in \Ob(\Eng_s)$, the relation $\simeq^{N, N'}_s$ is an \textbf{equivalence relation} on the homset $\Hom_{\Eng_s}(N,N')$. That is, for any aspects $V,V',V'' \in \Hom_{\Eng_s}(N, N')$
\begin{itemize}[topsep=-1ex,itemsep=-1.5ex]
\item[\textbf{(a)}] $\simeq_s^{N,N'}$ is reflexive: $V \simeq^{N, N'}_s V$ . 

\item[\textbf{(b)}] $\simeq_s^{N, N'}$ is symmetric: $V \simeq^{N, N'}_s V'$ if, and only if, $V' \simeq^{N, N'}_s V'$.

\item[\textbf{(c)}] $\simeq_s^{N, N'}$ is transitive: if $V \simeq^{N, N'}_s V'$ and $V' \simeq^{N, N'}_s V''$, then $V \simeq^{N, N'}_s V''$. 
\end{itemize}
\item The collection $\simeq_s$ of equivalence relations $\simeq^{N, N'}_s$ defines a \textbf{congruence} on $\Eng_s$:
\begin{itemize}[topsep=-1ex,itemsep=-1.5ex]
\item[\textbf{(d)}] For every $N, N', N''$, if $V_1 \simeq^{N, N'}_s V'_1$ and $V_2 \simeq^{N', N''}_s V'_2$, then $V_1; V_2 \simeq^{N, N''}_s V'_1; V'_2$. 
\[ 
\begin{tikzpicture}[description/.style={fill=white,inner sep=2pt}]
\matrix (m) [matrix of math nodes, row sep=1.15em, column sep=1.5em, text height=1.75ex, text depth=1.25ex]
{ \N & \checkmark & \N' & \checkmark & \N'' \\ };
\path[->]
(m-1-1) edge [bend left=30] node[above] {\footnotesize$\V_1$} (m-1-3) edge [bend right=30] node[below] {\footnotesize$\V'_1$} (m-1-3)
(m-1-3) edge [bend left=30] node[above] {\footnotesize$\V_2$} (m-1-5) edge [bend right=30] node[below] {\footnotesize$\V'_2$} (m-1-5);
\end{tikzpicture} 
\]

\vspace{-1cm} \mbox{}
~%
\lozend
\end{itemize}
\end{itemize}
\end{proposition}

\begin{remark} 
Although Propositions~\ref{prop:Engs} and \ref{prop:congruence} are stated in terms of a single author $s$, they are also valid for any fixed author set $S$. ~%
\lozend
\end{remark}


\subsubsection*{Proof that $\simeq_s$ is a congruence relation on $\Eng_s$ [technical]} 

Proposition~\ref{prop:congruence} is a consequence of the linguistic postulates of \S~\ref{sec:endoresement_postulates}. For \textbf{(a)}, if $V \colon N \to N'$ is an aspect such that $s$ endorses $\V$ as functional, then the syntax formula in Figure~\ref{fig:FormulaFact} is true:
\begin{quote}
For any $\qt{\N}$ $x$, we know that $x$ $\qt{\V \N'}$, that we call $y_1$, and we know that $x$ $\qt{\V \N'}$, that we call $y_2$; and the fact is, $y_1$ and $y_2$ are the same for any $x$.
\end{quote}
Note that the functionality of $\V$ implies that $y_1 = y_2$. Condition \textbf{(b)}, on the other hand, follows directly by applying the same syntax formula, since ``$=$" is a symmetric relation. Condition \textbf{(c)} follows similarly. 

Condition \textbf{(d)} is more delicate to check. Suppose we are given aspects $V_1, V'_1 \colon N \to N'$ and $V_2, V'_2 \colon N' \to N''$ such that $s \models [V_1 \simeq V'_1]$ and $s \models [V_2 \simeq V'_2]$. Let $x$ be any $\qt{\N}$. Then then following statement is true:
\begin{quote}
We know that $x \qt{\V_1 \N'}$, that we call $y_1$, and we know that $x \qt{V'_1 \N'}$, that we call $y_2$; and the fact is, $y_1$ and $y_2$ are the same.
\end{quote}
Set $y = y_1 = y_2$. Since $s \models [V_2 \simeq V'_2]$, we have that the following statement is true for $y$:
\begin{quote}
We know that $y \qt{\V_2 N''}$, that we call $z_1$, and we know that $y \qt{V'_2 N''}$, that we call $z_2$; and the fact is, $z_1$ and $z_2$ are the same. 
\end{quote}

Now, consider the following statement
\begin{quote}
We know that for any $\qt{\N}$ $x$, we know that $x \qt{\V_1 \N'}$, which $\qt{\V_2 \N''}$, that we call $w_1$, and we know that $x \qt{\V'_1 \N'}$, which $\qt{\V'_2 \N''}$, that we call $w_2$; and the fact is, $w_1$ and $w_2$ are the same.
\end{quote}
Since $s \models \V_1; \V_2$ and $s \models \V'_1; \V'_2$ by Postulate~\ref{rules:EndorsementConcatenation}, we have that $w_1 = z_1$ and $w_2 = z_2$. On the other hand, we know that $z_1 = z_2$. It follows that the previous statement is true, and hence $s \models [\V_1; \V_2 \simeq \V'_1; \V'_2]$.


\subsubsection*{Structure of $\Eng$ as a bicategory [technical]}

Denote by $\Ob(\Eng)$ the collection of types in $\Eng$. In order to see that \textbf{(a)}, \textbf{(b)} and \textbf{(c)} from Definition~\ref{def:type_aspect_fact} define a bicategory $\Eng$, the first thing we need to know is that for each pair of types $N_1, N_2 \in \Ob(\Eng)$, there is a category $\Eng(N_1, N_2)$, whose objects are given by aspects $V \colon N_1 \to N_2$ (\ie our $1$-cells), and whose morphisms between aspects $V_1, V_2 \colon N_1 \to N_2$ (\ie our $2$-cells) are given by the facts $F \colon V_1 \Rightarrow V_2$.  
\begin{figure}[H]
\centering
\belowcaptionskip = -15pt
\begin{tikzpicture}[description/.style={fill=white,inner sep=2pt}]
\matrix (m) [matrix of math nodes, row sep=0.25em, column sep=0.75em]
{ N_1 & \Downarrow F & N_2 \\ };
\path[->]
(m-1-1) edge [bend left=45] node[above] {\footnotesize$V_1$} (m-1-3)
(m-1-1) edge [bend right=45] node[below] {\footnotesize$V_2$} (m-1-3);
\end{tikzpicture}
\caption{$2$-cells in $\Eng$.}\label{fig:2_cells_Eng}
\end{figure}
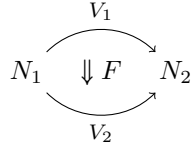
The composition of facts in $\Eng(N_1, N_2)$ is called \textbf{vertical composition}: given two facts $F_1 \colon V_1 \Rightarrow V_2$ and $F_2 \colon V_2 \Rightarrow V_3$, the composition is given by intersecting author sets:
\begin{equation}\label{eqn:Eng_vertical}
F_2 \circ F_1 := F_1 \cap F_2 \models [\V_1 \simeq \V_3].
\end{equation}
This composition defines a fact by Proposition~\ref{prop:congruence} \textbf{(c)}, and is associative since the intersection of sets is. Now, for every aspect $V \colon N_1 \to N_2$ there is a unique identity fact at $V$, given by the endorsement $\Auth(V) \models [(\N_1, \V, \N_2) \simeq (\N_1, \V, \N_2)]$ (see Proposition~\ref{prop:congruence} \textbf{(a)}). It follows by definition that these facts, that we denote ${\rm id}_{V}$, are unital with respect to the previous composition. Therefore, $\Eng(N_1, N_2)$ defines a category.

To define \textbf{horizontal composition} we will need to supply, for every $N_1, N_2, N_3 \in \Ob(\Eng)$, a functor 
\[ 
\Eng(N_1, N_2) \times \Eng(N_2, N_3) \too \Eng(N_1, N_3).
\]
We will denote it by
\begin{align*}
\big( V, V' \big) & \mapsto V; V' \;\;\mbox{ on aspects}, \\
(F, F') & \mapsto F \bullet F' \;\;\mbox{ on facts}.
\end{align*}
Horizontal composition of aspects $V=(\V,\Auth(V))$ works by intersecting author sets and concatenating verb phrases, as in Definition~\ref{def:Concatenation}, \ie
\begin{equation}\label{eqn:Eng_horizontal_aspects}
V; V'= \big((\V; \V'), \Auth(V) \cap \Auth(V)\big)
\end{equation}
Let $\Auth(V;V')\coloneqq\Auth(V)\cap\Auth(V')$. 

Vertical composition of facts is given as follows. If we are given $F \colon V_1 \Rightarrow V_2$ and $F' \colon V'_1 \Rightarrow V'_2$, define a fact $F \bullet F' \colon V_1;V'_1 \Rightarrow V_2;V'_2$ to be the endorsed equivalence
\begin{equation}\label{eqn:Eng_horizontal_facts}
F\cap F'\models \big[ (\V_1; \V_2) \simeq (\V'_1; \V'_2) \big]
\end{equation}
(see Proposition~\ref{prop:congruence} \textbf{(d)}). Note that this composition is associative at the level of aspects (see Remark~\ref{rem:concatenation_associative}). Moreover, for each type $N \in \Ob(\Eng)$, there is an identity aspect ${\rm id}_{N}$ in $\Eng(N, N)$ given by the verb phrase $e$, as in Postulate~\ref{rules:I}, and endorsed by $\Auth(N)$. Note that ${\rm id}_{N}$ is unital with respect to the horizontal composition at the level of aspects. At the level of facts, we can check that for each aspect $V \colon N_1 \to N_2$, there are isomorphisms 
\[
V; {\rm id}_{N_2} \Rightarrow V \mbox{ \ and \ } {\rm id}_{N_1}; V \Rightarrow V, 
\] 
since the equivalences $[\V; e] \simeq \V$ and $[e; \V] \simeq \V$ are endorsed by $\Auth(V) \cap \Auth(N_2) = \Auth(V)$ and $\Auth(V) \cap \Auth(N_1) = \Auth(V)$ by Postulate~\ref{rules:I_unital}. 

Finally, according to \cite[Definitions 1.2.5 \& 1.5.1]{Leinster}, it is only remains to verify that the previous associativity and unital isomorphisms satisfy the pentagon and triangle axioms (coherence conditions) to conclude that $\Eng$ is a bicategory, but this is a straightforward consequence of the previous constructions. Summarizing, we have proved Proposition~\ref{prop:Eng}.

\begin{proposition}\label{prop:Eng} 
The collection $\Eng$ of all types, aspects, and facts as in Definition~\ref{def:type_aspect_fact}, equipped with the vertical and horizontal compositions given by (\ref{eqn:Eng_vertical}), (\ref{eqn:Eng_horizontal_aspects}), and (\ref{eqn:Eng_horizontal_facts}), forms a bicategory. ~%
\lozend
\end{proposition}

\begin{remark}\label{rem:dual}
Let $V_1 \Rightarrow V_2$ be a $2$-cell in $\Eng$, \ie a fact. The direction ``$\Rightarrow$" is determined by the way the equivalence $\V_1 \simeq \V_2$ is read according to Figure~\ref{fig:FormulaFact}, \ie we read the statement involving ``\qt{\V_1}" first. Note by the symmetry of the relation $\simeq$ (see Proposition~\ref{prop:congruence}), that if $V_1\Rightarrow V_2$ is a fact with author set $F$, then so is $V_2\Rightarrow V_1$. ~%
\lozend
\end{remark}


\subsection{The bicategory $\iEng$ of instantiated English expressions}\label{sec:Eng_instances}

In this section we add to $\Eng$ some additional structure by exemplifying each linguistic expression in it. That is, we bundle each type $N$ with a set of examples, or as we call them, tokens,\footnote{The motivation behind the term ``token" comes from the usual usage in philosophy, \eg \cite{TyTo}. Tokens are real-world examples and instances of abstract types. In this paper the word ``example" has a wider connotation, while the term ``instance" will be reserved for \Defcolor~\ref{def:Instantiation}.} where the tokens must match the noun phrase $\qt{\N}$. Thus, we obtain a new bicategory $\iEng$ of instantiated English expressions, as shown in Proposition~\ref{prop:iEng}. Its $0$-cells, $1$-cells and $2$-cells will be described in Definition~\ref{def:type_aspect_fact_instance}. The bicategory $\iEng$ will provide the necessary theoretical framework to present the concept of instantiated olog. 

Just as for linguistic expressions, we need some postulates about author endorsement to govern the behavior of tokens with respect to noun phrases and verb phrases. For example, if you understand the noun phrase \fakebox{a US president}, then you probably endorse that Abraham Lincoln is an example, or token, of the corresponding set.

Recall the notion of endorsement for noun phrases and functional verb phrases from Definition~\ref{def:EndorsementNoun}.

\begin{definition}\label{def:validation} \
\begin{itemize}[topsep=-1ex,itemsep=-1.5ex]
\item[\textbf{(a)}] Let $\N$ be a noun phrase. One says that $x$ is a \textbf{token} of $\N$ if $x$ is an example to which the noun phrase $\qt{\N}$ applies, \ie if the following sentence is true:
\begin{equation}\label{eqn:endorsement_token_noun}
\N(x):=\;\; x\tn{  ``is" }\qt{\N}.
\end{equation}
If $s$ is an endorsing author $s \models \N$ and $x$ is a token of $\N$ according to $s$, we say that $s$ \textbf{endorses $x$ as a token of $\N$}; this will be denoted by $s \models (x : \N)$.

\item[\textbf{(b)}] Let $(\N_1, \V, \N_2)$ be a sentence, and let $x$ and $y$ be tokens of $\N_1$ and $\N_2$, respectively. One says that \textbf{$x$ corresponds to $y$ via $\V$} if the following sentence, denoted $\V(x,y)$, is true:
\begin{equation}\label{eqn:correspondence}
\V(x,y) := \;\;x \tn{ ``is" } \qt{\N_1} \tn{ ``which" } \qt{\V \N_2}, \tn{ ``namely'' } y.
\end{equation}
If $s$ is endorsing author $s \models (\N_1, \V, \N_2), s \models (x : \N_1), s \models (y : \N_2),$ and if $x$ corresponds to $y$ via $\V$ according to $s$, then we say that $s$ \textbf{endorses the correspondence $\V(x, y)$ between $x$ and $y$}; this will be denoted by $s \models \V(x,y)$.
\end{itemize}
If $S$ is an author set, we say that $S$ endorses either a token $S\models (x:\N)$ or a correspondence between tokens $S\models\V(x,y)$, if every member of $s$ endorses it, \ie $s\models (x:\N)$ or $s\models \V(x,y)$ for every $s\in S$.
\lozend
\end{definition}

\begin{example}
Consider the following sentence, which we endorse as functional:
$$\ovbox{P}{a US president}\To{\parbox{2in}{\small\centering can, using the usual chronological ordering of presidents, be assigned}}\ovbox{I}{an integer}$$
We endorse that Abraham Lincoln is a token of $P$, that $16$ is a token of $I$, and that they correspond via the verb phrase labeling the arrow. That is, we agree that: ``Abraham Lincoln is a US president which can, using the usual chronological ordering of presidents, be assigned an integer, namely 16."
~\lozend
\end{example}

We now postulate that the correspondence (\ref{eqn:correspondence}) is functional, and that the unit verb phrase and composition of verb phrases act like the identity function and function composition for instances, respectively.

\begin{postulate}\label{rules:EndorsementFact} 
Let $s$ be a person who understands Definition~\ref{def:validation}.
\begin{itemize}[topsep=-1ex,itemsep=-1.5ex]
\item[\textbf{(a)}] Suppose given an endorsement $s \models (\N_1, \V, \N_2)$. Then for every token $x$ of $N_1$ such that $s \models (x : N_1)$, there exists a unique token $y$ of $N_2$ such that $s \models \V(x,y)$. 

\item[\textbf{(b)}] Suppose given endorsements for concatenatable sentences, $s \models (\N_1, \V_1, \N_2)$ and $s \models (\N_2, \V_2, \N_3)$. Note that  $s \models (\N_1, \V_1; \V_2, \N_3)$ by Postulate~\ref{rules:EndorsementConcatenation}. Suppose also that $s$ endorses tokens $s \models (x : \N_1), (y : \N_2), (z : \N_3)$ and correspondences $s \models \V_1(x,y)$ and $s \models \V_2(y,z)$. Then we will assume that $s \models (\V_1; \V_2)(x,z)$. 

\item[\textbf{(c)}] Suppose given an endorsement $s \models (x : \N)$ for some noun phrase $\N$. Then we will assume that $s \models \unit(x, x)$, with $\unit$ as the unit verb phrase defined in Postulate~\ref{rules:I}. \lozend
\end{itemize}
\end{postulate}

Recall the definitions of types, aspects, and facts from Definition~\ref{def:type_aspect_fact}, and the notion of tokens and their correspondences from Definition~\ref{def:validation}.

\begin{definition}\label{def:type_aspect_fact_instance} \
\begin{itemize}[topsep=-1ex,itemsep=-1.5ex]
\item[\textbf{(a)}] An \textbf{instantiated type} is a pair $\oN=(N,\T)$, where $N$ is a type with noun phrase $\N$ and $\T$ is a set  such that for every $x \in \T$, we have an endorsement 
$$\Auth(N) \models (x : \N)$$ 
that $x$ is a token of $\N$. We may write $\T(\oN)\coloneqq \T$ and refer to it the \textbf{token set} of $\oN$.

\item[\textbf{(b)}] Suppose given two instantiated types $\oN_1$ and $\oN_2$. An \textbf{instantiated aspect} from $\oN_1$ to $\oN_2$, denoted $\oV \colon \oN_1 \to \oN_2$, consists of a pair $(V,f)$, where $V\colon N_1\to N_2$ is an aspect, and $f$ is a function $f \colon \T(\oN_1) \to \T(\oN_2)$, such that for every token $x \in \T(\oN_1)$, there is an endorsement 
$$\Auth(V)\models \V(x, f(x)),$$ 
where $\V$ is the verb phrase of $V$. We may write $\T(\oV)\coloneqq f$ and refer to it as the \textbf{token function} of $\oV$. 

\item[\textbf{(c)}] Suppose given two instantiated aspects $\oV_1, \oV_2 \colon \oN_1 \to \oN_2$ with equal token functions, $\T(\oV_1) = \T(\oV_2)$. An \textbf{instantiated fact} from $\oV_1$ to $\oV_2$, denoted $\oF\colon\oV_1 \Rightarrow \oV_2$, consists of a fact $\oF=F$, where
$$F \models [\V_1 \simeq \V_2]$$
is a fact in $\Eng$.  ~%
\lozend
\end{itemize}
\end{definition}

As in \S~\ref{sec:Eng}, one can prove that the collections of instantiated types, aspects and facts in Definition~\ref{def:type_aspect_fact_instance} form a bicategory $\iEng\ss\Eng\times\Set$. We do not give many details on this, since the proof proceeds as it did with $\Eng$ in the previous section. However, we will specify how vertical and horizontal compositions are defined.  On a first reading, one may skip to Proposition~\ref{prop:iEng}.


\subsubsection*{Structure of the bicategory $\iEng$ [technical]}

Denote by $\Ob(\iEng)$ the collection of instantiated types, as in Definition~\ref{def:type_aspect_fact_instance}. For each pair $\oN_1, \oN_2 \in \Ob(\iEng)$, there is a category $\iEng\big(\oN_1, \oN_2 \big)$ whose objects (\ie $1$-cells in $\iEng$) are given by instantiated aspects $\oV \colon \oN_1 \to \oN_2$, and whose morphisms (\ie $2$-cells in $\iEng$) $\oV_1 \Rightarrow \oV_2$ are given by:
\[
\Hom_{\iEng\big( \oN_1, \oN_2 \big)}\big( \oV_1, \oV_2 \big) :=
\left\{ 
\begin{array}{ll} 
\Hom_{\Eng\big( N_1, N_2 \big)}\big( V_1, V_2 \big) & \mbox{if $\T(\oV_1) = \T(\oV_2)$}, \\
\emptyset & \mbox{otherwise}. 
\end{array} 
\right.
\]
By the previous formula it follows that the \textbf{vertical composition} in $\iEng$ is just the vertical composition in $\Eng$, and identities are also straightforward.

For every triple $\oN_1, \oN_2, \oN_3 \in \Ob(\iEng)$, there is a \textbf{horizontal composition} functor 
\[
\iEng(\oN_1, \oN_2) \times \iEng(\oN_2, \oN_3) \to \iEng(\oN_1, \oN_3),
\]
denoted $\oV_1; \oV_2$ for every $\oV_1 \in \Ob(\iEng(\oN_1, \oN_2))$ and $\oV_2 \in \Ob(\iEng(\oN_2, \oN_3))$, acting on instantiated aspects according to the intersection of sets, concatenation of verb phrases, and composition of functions.

At the level of instantiated facts, if we are given $\oF \colon \oV_1 \Rightarrow \oV_2$ and $\oF' \colon \oV'_1 \Rightarrow \oV'_2$, then the instantiated fact
\[
\oF \bullet \oF' \colon \oV_1; \oV_2 \Rightarrow \oV'_1; \oV'_2,
\]
is given by the endorsed equivalence $\oF \cap \oF' \models [\oV_1; \oV'_1 \simeq \oV_2; \oV'_2]$. Finally, for each instantiated type $\oN$, there is an identity instantiated aspect ${\rm id}_{\oN}$ given by unit the noun phrase $e$ as in Postulate~\ref{rules:I}, endorsed by $\Auth(N)$, and by the identity function ${\rm id}_{\T(\oN)} \colon \T(\oN) \to \T(\oN)$, which is unital with respect to the horizontal composition of instantiated aspects. Note that the previous operations are well defined by Postulate~\ref{rules:EndorsementFact}.

\begin{proposition}\label{prop:iEng} 
The collection of all types, aspects, and facts as in Definition~\ref{def:type_aspect_fact_instance} forms a bicategory, denoted $\iEng$.~%
\lozend
\end{proposition}

As was true with $\Eng$, if we fix a single author $s$ (resp. an author set $S$), there is a subcategory $\iEng_s\ss\iEng$ of instantiated linguistic expressions endorsed by $s$ (resp. $S$). There is no need to consider 2-cells in $\iEng_s$, and composition is given simply by concatenating verb phrases and composing functions. 

The following result is just a matter of checking the conditions in \cite[Definition~1.5.8]{Leinster}.

\begin{proposition}\label{prop:Instantiation} \
\begin{itemize}[topsep=-1ex,itemsep=-1.5ex]
\item[\textbf{(a)}] The mappings
\begin{align*}
\oN & \mapsto N, \mbox{ for every instantiated type $\oN$, and} \\
\oV & \mapsto V, \mbox{ for every instantiated aspect $\oV$}.
\end{align*}
define an strict functor between bicategories $\forget \colon \iEng \too \Eng$. 

\item[\textbf{(b)}] The mappings
\begin{align*}
\oN & \mapsto \T(\oN), \mbox{ for every instantiated type $\oN$, and} \\
\oV & \mapsto \T(\oV), \mbox{ for every instantiated aspect $\oV$},
\end{align*}
define an strict functor between bicategories $\T \colon \iEng \too \Set$.
\end{itemize}
Together, these define an inclusion $(U,\T)\colon\iEng\hookrightarrow\Eng\times\Set$. ~%
\lozend
\end{proposition}


\section{Ologs, instantiated ologs, and mappings between them}\label{sec:def}

We begin this section formalizing the definition of ologs found in \cite{SpivakOlogs}, using a new construction, which we will call a {\it linguistic structure} in Definition~\ref{def:Linguistic}. We suggest that readers who are unfamiliar with ologs consult \cite[\S 1, 2, 3, and 4]{SpivakOlogsPre} or \cite[\S 1, 2, and 3]{SpivakOlogs}.%
\footnote{Note that the later sections of \cite{SpivakOlogsPre}, \cite{SpivakOlogs} also discuss pullbacks and pushouts within an olog; we do not address these notions in this paper.}

In \S~\ref{sec:com}, we define mappings between ologs, called \textit{linguistic functors} (see Definition~\ref{def:LinguisticFunctor}). These are mappings between categories that preserve the linguistic structure of the ologs. In Example \ref{ex:no_chance}, we explain this by providing two functors $F$ and $G$, between the same two ologs, such that $F$ is linguistic and $G$ fails to be. In \S~\ref{sec:instantiations}  we define an instantiated olog to be a standard olog equipped with a set of examples that conforms to its linguistic structure. This formalizes the original notion of olog given in \cite{SpivakOlogs}. Finally, in \S~\ref{sec:instantiated_functors} we study mappings between instantiated ologs, called \textit{instantiated functors}, which are instantiated analogues of linguistic functors. In particular, we provide a more linguistic version of the notion of \textit{meaningful functor} given in \cite{SpivakOlogs}. All of the above is achieved by appealing to our definitions of $\Eng$ and $\iEng$ from \S~\ref{sec:EngBicat}.


\subsection{Linguistic structures}\label{sec:ling} 

Recall that an olog, as defined in \cite{SpivakOlogs}, is a category whose objects and morphisms are labeled with {\it noun phrases} and {\it verb phrases}, respectively, in order to model a conceptual situation. Commutative diagrams in an olog are called \textit{facts}\footnote{Notice that in our proposed definition of facts, equivalences between sentences must be endorsed by an author set, and this is not explicitly mentioned in \cite{SpivakOlogs}.}; they are equivalences between two sentences in the English language. Importantly, these types, aspects, and facts must follow certain guidelines, to ensure that the category-theoretic meaning is aligned with the conceptual and linguistic meaning intended by the authors. The {\it Rules of Good Practice} stated in \cite[2.1.2, 2.2.3, and 2.3.4]{SpivakOlogs} are such a set of guidelines. An example of these rules is that the label of each arrow $a \rightarrow b$ corresponds to a mathematical function \cite{WikiFun}. All of these notions are author-dependent, as we made explicit using the notion of \emph{endorsement}, as in Definition~\ref{def:EndorsementNoun}. After our work defining $\Eng$ in \S~\ref{sec:EngBicat}, formalizing the notion of olog given in \cite{SpivakOlogs} is straightforward.

\begin{definition}\label{def:Linguistic} 
Let $\cat{C}$ be a category. A \textbf{linguistic structure on $\cat{C}$} is a lax functor $L \colon \cat{C} \too \Eng$. An \textbf{olog} is a pair $(\cat{C}, L)$, where $\cat{C}$ is a category, and $L$ is a linguistic structure on $\cat{C}$. ~%
\lozend
\end{definition}

The notion of lax functor, found in Definition~\ref{def:Linguistic}, is standard in category theory literature; see for example \cite[Definition 1.5.8]{Leinster}. However, we will now spell it out explicitly in our case.

\begin{remark}\label{rem:Linguistic} 
Let $L \colon \cat{C} \too \Eng$ be a linguistic structure on $\cat{C}$. Then: 
\begin{itemize}[topsep=-1ex,itemsep=-1.5ex]
\item[\textbf{(a)}] Every $c \in \Ob(\cat{C})$ is mapped into a type $L(c)$, consisting of a noun phrase $\mfL(c)$ and an endorsement $\Auth(L(c)) \models \mfL(c)$.

\item[\textbf{(b)}] Every $f \in \Hom_{\cat{C}}(c,c')$ is mapped into an aspect $L(f)$, consisting of a functional verb phrase $\mfL(f) \colon \mfL(c) \to \mfL(c')$ and an endorsement $\Auth(L(f)) \models \mfL(f)$.

\item[\textbf{(c)}]  Let $a\To{f}b\To{g}c$ be a pair of composable morphisms in $\cat{C}$, with $h=f;g$. Then $L$ being lax means that there 2-cell $F \colon L(h) \Rightarrow L(f); L(g)$, \ie a fact, 
\[
F \models [ \mfL(h) \simeq \mfL(f); \mfL(g)].
\]
\end{itemize}
\vspace{-1.2cm}\lozend
\end{remark}

\begin{example}\label{ex:olog} 
We run through \Defcolor~\ref{def:Linguistic} in the case of a category $\cat{C}$,
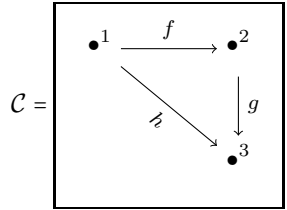
\[ 
\cat{C} = \parbox{1.25in}{
\fbox{
\begin{tikzpicture}[description/.style={fill=white,inner sep=2pt}]
\matrix (m) [ampersand replacement=\&, matrix of math nodes, row sep=0.1em, column sep=1.5em, text height=1.75ex, text depth=1.25ex]
{ \bullet^1 \& \& \bullet^{2} \\ \& \mbox{} \\ \& \& \bullet^{3} \\ };
\path[->]
(m-1-1) edge node[below,sloped] {\footnotesize$h$} (m-3-3)
(m-1-1) edge node[above] {\footnotesize$f$} (m-1-3)
(m-1-3) edge node[right] {\footnotesize$g$} (m-3-3);
\end{tikzpicture} 
}}
\]
and a linguistic structure $L$ on it, depicted here:
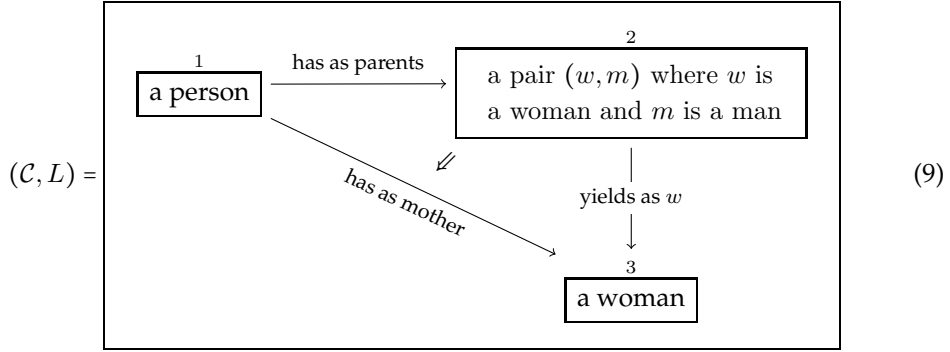
\begin{align}\label{dia:linguistic structure objects explicit}
(\cat{C}, L) = \parbox{3.8in}{
\fbox{ \begin{tikzpicture}[description/.style={fill=white,inner sep=2pt}] 
\matrix (m) [ampersand replacement=\&, matrix of math nodes, row sep=1.6em, column sep=3em] 
{ \ovbox{1}{a person} \& \& \ovbox{2}{ $\begin{array}{ll} \mbox{{\rm a pair $(w,m)$ where $w$ is}} \\ \mbox{{\rm  a woman and $m$ is a man}} \end{array}$} \\ \& \mbox{} \\ \mbox{} \& \& \ovbox{3}{a woman} \\ }; 
\path[->] 
(m-1-1) edge node[above] {$\footnotesize\mbox{has as parents}$} (m-1-3) edge node[below,sloped] {$\footnotesize\mbox{has as mother}$} (m-3-3)
(m-1-3) edge node[description] {$\footnotesize\mbox{yields as $w$}$}  (m-3-3)
(m-1-3)-- node[pos=0.3] {$\Swarrow$} (m-2-2); 
\end{tikzpicture} 
}}
\end{align}
Each object has been assigned a noun phrase, and each arrow has been assigned a verb phrase. The two paths from $1$ to $3$, namely $h$ and the composition of $f$ followed by $g$, have been declared equivalent. 

We have not yet specified the author sets for the types, aspects, and facts above. This is arbitrary, but we specify a full example for completeness. Suppose we say 
$$\Auth(1)\coloneqq\{A,B,C\},\quad \Auth(2)\coloneqq\{A\},\quad \Auth(3)\coloneqq\{A,B\}$$
Then we must have $\Auth(f)=\emptyset$ or $\Auth(f)=\{A\}$ because it must be a subset, $\Auth(f)\ss\Auth(1)\cap\Auth(2)$; similarly for $\Auth(g)$. Let's say 
$$\Auth(f)=\Auth(g)\coloneqq\{A\}\;\;\tn{ and }\;\;\Auth(h)=\{B\}.$$ 
Since $\Auth(f;g)=\Auth(f)\cap\Auth(g)=\{A\}$, but $\Auth(h)=\{B\}$, we must have $F=\emptyset$. That is, given our choices, there is no one who is capable of endorsing the equivalence $\mfL(h) \simeq \mfL(f); \mfL(g)$. However, it is still a fact because the diagram commutes in the ground category $\cat{C}$. ~%
\lozend
\end{example}


\subsection{Mappings between ologs}\label{sec:com}

In this section we define the notion of morphism between ologs, which we call a \emph{linguistic functor}, in Definition~\ref{def:LinguisticFunctor}. After this, we can define the category of ologs, denoted $\Olog$. We then give in Example~\ref{ex:no_chance} two ologs and two functors between them, one of which makes linguistic sense and the other does not, which recapitulates the issue we presented in our introductory example (\ref{eqn:no_chance}).

\begin{definition}\label{def:LinguisticFunctor} 
Let $(\cat{C}, L)$ and $(\cat{D}, M)$ be ologs, as in Definition~\ref{def:Linguistic}. A \textbf{linguistic functor} between them, denoted $(F, F^{\#}) \colon (\cat{C}, L) \too (\cat{D}, M)$ consists of a functor $F \colon \cat{C} \too \cat{D}$ together with a lax transformation $F^{\#} \colon L \leadsto M \circ F$ as shown below:
\begin{figure}[H]
\belowcaptionskip = -15pt
\centering
\begin{tikzpicture}[description/.style={fill=white,inner sep=2pt}]
\matrix (m) [matrix of math nodes, row sep=0.15em, column sep=2.75em, text height=1.75ex, text depth=1.25ex]
{ \cat{C} & & \cat{D} \\ & \overset{F^{\#}}{\leadsto} \\ & \Eng \\ };
\path[->]
(m-1-1) edge node[above] {\footnotesize$F$} (m-1-3) edge node[left] {\footnotesize$L$} (m-3-2)
(m-1-3) edge node[right] {\footnotesize$M$} (m-3-2);
\end{tikzpicture}
\caption{Linguistic functor.}\label{fig:linguistic_functor}
\end{figure}
If $F^{\#}$ is the identity, in which case $L=M\circ F$, we say that the linguistic structure $L$ is the \textbf{pullback of $M$ along $F$} and write $L=F^*(M)$.

We define the \textbf{category of ologs}, denoted $\Olog$, to be the category whose objects are ologs and whose morphisms are linguistic functors (see Remark~\ref{rem:Grothendieck}). \lozend
\end{definition}

\begin{remark} \label{rem:authorship_functor}
The morphism $F^{\#}$ includes, for each $c \in \Ob(\cat{C})$, an aspect $F^{\#}_c \colon L(c) \leadsto M(Fc)$, called the \textbf{$c$-component of $F^{\#}$}, with endorsing author set $\Auth(F^{\#}_c)$. And for every $f \in \Hom_{\mathcal{C}}(c,c')$, it includes a fact
\begin{equation}\label{eqn:morphism_ling_struc}
\parbox{2.1in}{
\begin{tikzpicture}[description/.style={fill=white,inner sep=2pt}] 
\matrix (m) [ampersand replacement=\&, matrix of math nodes, row sep=1em, column sep=1.25em]
{ L(c) \& \& M(Fc) \\ \& F^{\#}_f\Swarrow \\ L(c') \& \& M(Fc') \\ };
\path[->]
(m-1-1) edge [snake it] node[above] {\footnotesize$F^{\#}_c$} (m-1-3) edge node[left] {\footnotesize$L(f)$} (m-3-1)
(m-1-3) edge node[right] {\footnotesize$M(Ff)$} (m-3-3)
(m-3-1) edge [snake it] node[below] {\footnotesize$F^{\#}_{c'}$} (m-3-3);
\end{tikzpicture} 
}
\end{equation}
Both the arrows between objects in an olog, denoted $\to$, and the component arrows for maps between ologs, denoted $\leadsto$, are assigned functional verb phrases. Although we use differently-shaped arrows to denote them, the equivalence in~\eqref{eqn:morphism_ling_struc} is of the usual kind, as in Definition~\ref{def:DeclaredEq}. ~%
\lozend
\end{remark}

\begin{remark}\label{rem:Grothendieck} 

There is an evident functor $\Olog\to\Cat$, sending an olog $(\cat{C},L)$ to its underlying category $\cat{C}$. This is a Grothendieck fibration, meaning that if $F\colon\cat{B}\to\cat{C}$ is any functor, there is a notion of \textbf{cartesian morphism} of ologs over $F$. This is just the pullback along $F$, in the sense of Definition~\ref{def:LinguisticFunctor}. Grothendieck fibrations also come with a notion of vertical morphisms. These are the linguistic functors that restrict to the identity $\id_{\cat{C}}$ on the underlying category.\lozend
\end{remark}

\begin{example}\label{ex:mordos} 
Consider the following linguistic structures $L$ and $M$ on the same category, $\fbox{$\bullet^2 \longleftarrow \bullet^1 \longrightarrow \bullet^3$}$:

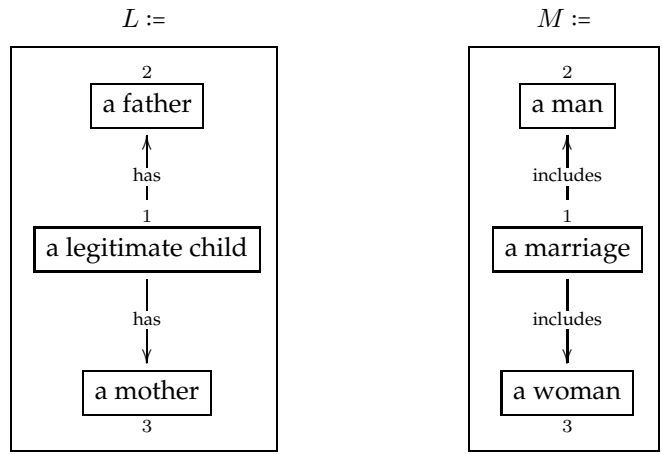
\begin{equation*}
\parbox{1.4in}{\boxtitle{$L :=$\\~}\fbox{
\xymatrix{
\ovbox{2}{a father}\\\ovbox{1}{a legitimate child}\ar[u]|{\tn{{\color{white}p}has{\color{white}p}}}\ar[d]|{\tn{{\color{white}p}has{\color{white}p}}}\\\unbox{3}{a mother}
}
}}
\hspace{1in}
\parbox{1in}{\boxtitle{$M :=$\\~}\fbox{
\xymatrix{
\ovbox{2}{a man}\\\ovbox{1}{a marriage}\ar[u]|{\tn{{\color{white}p}includes{\color{white}p}}}\ar[d]|{\tn{{\color{white}p}includes{\color{white}p}}}\\\unbox{3}{a woman}
}
}}
\end{equation*}

The phrase ``legitimate child" is an old-fashioned term for a child who was born in a marriage, which itself is required to be between a man and a woman. Let $S$ be an author set who endorses the linguistic expressions in $L$ and $M$. Suppose these authors also endorse the following sentences as functional:
\begin{center}
``a legitimate child was born in a marriage"\\
``a father is a man"\\
``a mother is a woman".
\end{center}
Technically speaking, they endorse the verb phrases \qt{\alpha_1}:= \mbox{``was born in"}; \qt{\alpha_2}:= \mbox{``is"}; and \qt{\alpha_3}:= \mbox{``is"} as components $\alpha_c$, for $c\in\{1,2,3\} = \Ob(\cat{C})$, of a lax transformation $\alpha \colon L \leadsto M$.

At this point, they have endorsed every object and arrow (both $\to$ and $\leadsto$) in the diagram below:
\begin{figure}[H]
\belowcaptionskip = -15pt
\centering
\begin{tikzpicture}[description/.style={fill=white,inner sep=2pt}]
\matrix (m) [matrix of math nodes, row sep=0.7em, column sep=1.75em]
{ \ovbox{2}{a father} & & \ovbox{2}{a man} \\
& \Nwarrow \\
\ovbox{1}{a legitimate child} & & \ovbox{1}{a marriage} \\
& \Swarrow \\
\unbox{3}{a mother} & & \unbox{3}{a woman} \\ };
\path[->]
(m-3-1) edge node[description] {$\footnotesize\mbox{has}$} (m-1-1) edge node[description] {$\footnotesize\mbox{has}$} (m-5-1) edge [snake it] node[above] {$\footnotesize \mbox{was born in}$} (m-3-3)
(m-1-1) edge [snake it] node[above] {$\footnotesize \mbox{is}$} (m-1-3)
(m-5-1) edge [snake it] node[below] {$\footnotesize \mbox{is}$} (m-5-3)
(m-3-3) edge node[description] {$\footnotesize\mbox{includes}$} (m-1-3) edge node[description] {$\footnotesize\mbox{includes}$} (m-5-3);
\end{tikzpicture}
\caption{Example of a morphism $\alpha\colon L\leadsto M$ between two linguistic structures on $\cat{C}$.}\label{fig:marriage!}
\end{figure}
The $2$-cells are drawn if each author $s\in S$ also endorses the corresponding equivalences. For example, a legitimate child has a father who is a man, and a legitimate child was born in a marriage, which includes a man. The point is, they had better be the same man(!) at least according to any authors endorsing $\alpha$. \lozend
\end{example}

\begin{remark}
Sometimes it is valuable to compare two linguistic structures $L$ and $M$ on the same category $\cat{C}$, even if one cannot satisfy the somewhat strong conditions of a morphism of linguistic structures. Different ologs on the same underlying category have been considered before in other fields, such as biology and materials science \cite{GiesaAnalogy}. In this case, one could consider a relaxed version of morphism $\alpha \colon L\leadsto M$ in which the component aspects $\alpha_c \colon L(c) \leadsto M(c)$ were not required to be functional, just ``conceptual", requiring a relaxed version of $\Eng$. We will not explain how to do this category-theoretically, although it is possible. Instead, we just give an example to clarify the idea.

Consider the left and right column, each a linguistic structure on the category \fbox{$\bullet\to\bullet$}. 
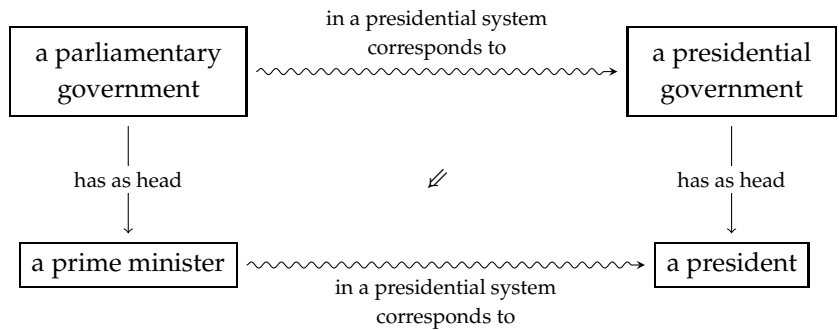
\[
\begin{tikzpicture}[description/.style={fill=white,inner sep=2pt}]
\matrix (m) [matrix of math nodes, row sep=1.25em, column sep=6em]
{ \obox{$\begin{array}{c} \mbox{a parliamentary} \\ \mbox{government} \end{array}$} & & \obox{$\begin{array}{c} \mbox{a presidential}  \\ \mbox{government} \end{array}$} \\ & \Swarrow \\ \obox{a prime minister} & & \obox{a president} \\ };
\path[->]
(m-1-1) edge [snake it] node[above] {\footnotesize$\begin{array}{c} \mbox{in a presidential system} \\ \mbox{corresponds to} \end{array}$} (m-1-3) edge node[description] {\footnotesize$\mbox{has as head}$} (m-3-1)
(m-1-3) edge node[description] {\footnotesize$\mbox{has as head}$} (m-3-3)
(m-3-1) edge [snake it] node[below] {\footnotesize$\begin{array}{c} \mbox{in a presidential system} \\ \mbox{corresponds to} \end{array}$} (m-3-3);
\end{tikzpicture}
\]
We cannot endorse the component maps between them (labeled, ``in a presidential system corresponds to") as being functional. However, the diagram still shows a valid analogy between two forms of government. Analogies have been considered from a category-theoretic perspective before; see \cite{Brown}. \lozend
\end{remark}

\begin{example}\label{ex:no_chance} 
In this example we will show two ologs $(\cat{C}, L)$ and $(\cat{D}, M)$, and two functors $F, G \colon \cat{C} \too \cat{D}$ between their underlying categories. We will find that it is straightforward to put a linguistic structure on $F$. On the other hand, we will explain why there is very little chance any author will endorse a linguistic structure on $G$. This way, we will explain an example from the introduction, \ssec~\ref{sec:intro}, namely the lack of linguistic structure on the functor in (\ref{eqn:no_chance}).

Consider the following ologs:
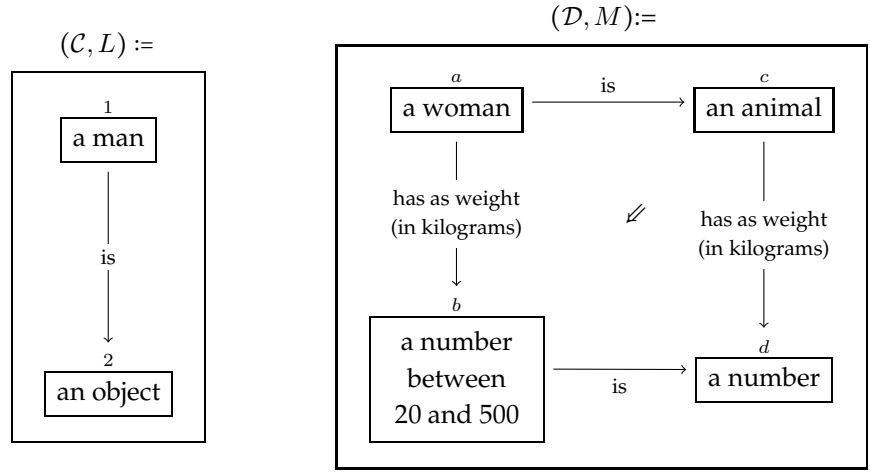
\begin{align*}
\parbox{1in}{\boxtitle{$(\cat{C}, L):=$\\~}\fbox{
\begin{tikzpicture}[description/.style={fill=white,inner sep=2pt}] 
\matrix (m) [ampersand replacement=\&, matrix of math nodes, row sep=6.5em]
{ \ovbox{1}{a man}\\ \ovbox{2}{an object} \\ }; 
\path[->] 
(m-1-1) edge node[description] {$\footnotesize \mbox{is}$} (m-2-1); 
\end{tikzpicture} }
}
&\hspace{.7in}
\parbox{2.8in}{\boxtitle{$(\cat{D}, M)$:=\\~}\fbox{ 
\begin{tikzpicture}[description/.style={fill=white,inner sep=2pt}] 
\matrix (m) [ampersand replacement=\&, matrix of math nodes, row sep=5.5em, column sep=5em] 
{ \ovbox{a}{a woman} \& \ovbox{c}{an animal} \\ \ovbox{b}{$\begin{array}{c} \mbox{a number} \\ \mbox{between} \\ \mbox{20 and 500} \end{array}$} \& \ovbox{d}{a number} \\ }; 
\path[->] 
(m-1-1) edge node[above] {$\footnotesize \mbox{is}$} (m-1-2) 
(m-1-1) edge node[description] {$\footnotesize \begin{array}{c} \mbox{has as weight} \\ \mbox{(in kilograms)} \end{array}$} (m-2-1)
(m-1-2) edge node[description] {$\footnotesize \begin{array}{c} \mbox{has as weight} \\ \mbox{(in kilograms)} \end{array}$} (m-2-2)
(m-2-1) edge node[below] {$\footnotesize \mbox{is}$} (m-2-2)
(m-1-2)-- node[description] {$\Swarrow$} (m-2-1); 
\end{tikzpicture} }
}
\end{align*}
Let $F\colon \cat{C} \too\cat{D}$ be the unique functor such that $F(1) = c$ and $F(2)= d$, and let $G\colon \cat{C} \too \cat{D}$ be the unique functor such that $G(1) = a$ and $G(2)= b$. Then the pulled back linguistic structures $F^*(M)$ and $G^*(M)$ on $\cat{C}$ are shown on the left and right below. 
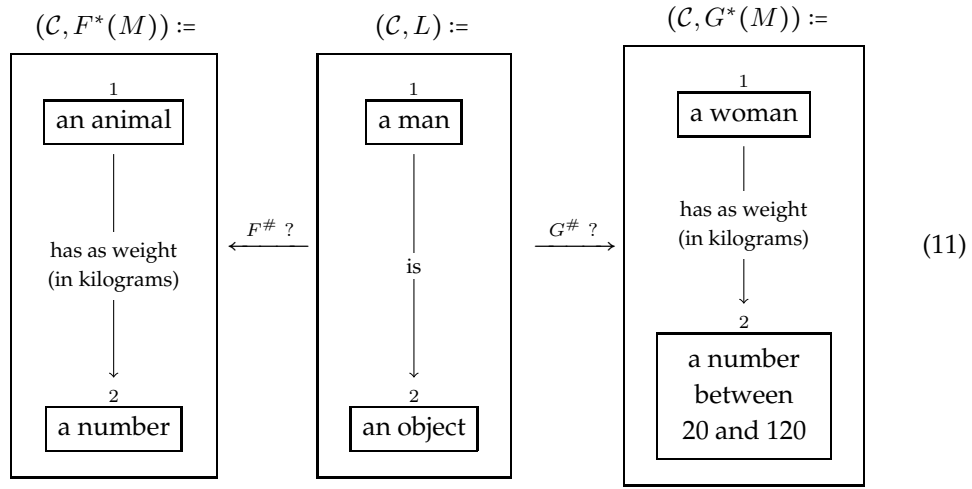
\begin{equation}\label{eqn:pullbacks_and_necessary}
\parbox{1.1in}{\boxtitle{$(\cat{C}, F^\ast(M)):=$\\~}\fbox{
\begin{tikzpicture}[description/.style={fill=white,inner sep=2pt}] 
\matrix (m) [ampersand replacement=\&, matrix of math nodes, row sep=8.5em]
{ \ovbox{1}{an animal}\\ \ovbox{2}{a number} \\ }; 
\path[->] 
(m-1-1) edge node[description] {$\footnotesize \begin{array}{c} \mbox{has as weight} \\ \mbox{(in kilograms)} \end{array}$} (m-2-1); 
\end{tikzpicture} }
}
\From{\;\; F^{\#}\ ?\;\;}
\parbox{1.1in}{\boxtitle{$(\cat{C}, L):=$\\~}\fbox{
\begin{tikzpicture}[description/.style={fill=white,inner sep=2pt}] 
\matrix (m) [ampersand replacement=\&, matrix of math nodes, row sep=8.5em]
{ \ovbox{1}{a man}\\ \ovbox{2}{an object} \\ }; 
\path[->] 
(m-1-1) edge node[description] {$\footnotesize \mbox{is}$} (m-2-1); 
\end{tikzpicture} }
}\To{\;\; G^{\#}\ ?\;\;}
\parbox{1.3in}{\boxtitle{$(\cat{C}, G^\ast(M)):=$\\~}\fbox{
\begin{tikzpicture}[description/.style={fill=white,inner sep=2pt}] 
\matrix (m) [ampersand replacement=\&, matrix of math nodes, row sep=6em]
{ \ovbox{1}{a woman}\\ \ovbox{2}{$\begin{array}{c} \mbox{a number} \\ \mbox{between} \\ \mbox{20 and 120} \end{array}$}\\ }; 
\path[->] 
(m-1-1) edge node[description] {$\footnotesize \begin{array}{c} \mbox{has as weight} \\ \mbox{(in kilograms)} \end{array}$} (m-2-1); 
\end{tikzpicture} }
}
\end{equation}
In order to extend $F$ and $G$ to linguistic functors, we need lax transformation $F^{\#} \colon L \leadsto M \circ F$ and $G^{\#} \colon L \leadsto M \circ G$, as indicated in (\ref{eqn:pullbacks_and_necessary}). We begin by endorsing a certain $F^{\#}$, and then return to explaining the difficulty of finding $G^{\#}$. 

To provide $F^{\#}$, we first need two component aspects; we use those labeling the top and bottom maps here:
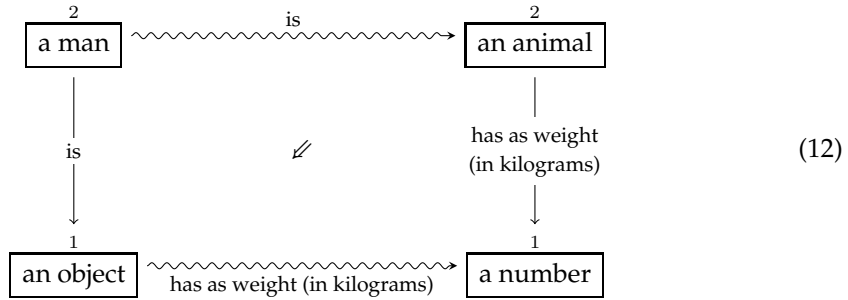
\begin{equation}\label{eqn:woman_number}
\parbox{3.4in}{
\begin{tikzpicture}[description/.style={fill=white,inner sep=2pt}]
\matrix (m) [ampersand replacement=\&, matrix of math nodes, row sep=2em, column sep=5em]
{ \ovbox{2}{a man} \& \& \ovbox{2}{an animal} \\
\& \Swarrow \\
\ovbox{1}{an object} \& \& \ovbox{1}{a number} \\ };
\path[->]
(m-1-1) edge node[description] {$\footnotesize\mbox{is}$} (m-3-1) edge [snake it] node[above] {$\footnotesize\mbox{is}  $} (m-1-3)
(m-1-3) edge node[description] {$\footnotesize\begin{array}{c} \mbox{has as weight} \\ \mbox{(in kilograms)} \end{array}$} (m-3-3) 
(m-3-1) edge [snake it] node[below] {$\footnotesize \mbox{has as weight (in kilograms)}$} (m-3-3);
\end{tikzpicture}
}
\end{equation}
It just suffices to ask whether, if we measure the weight (in kilograms) of a man, regarding him either as an object or as an animal, we get the same number. We endorse that fact, thus providing a linguistic functor $(F, F^{\#}) \colon (\cat{C}, L) \too (\cat{D}, M)$.

On the other hand, it is difficult to do the same for $G^{\#}$ as in the right side of \eqref{eqn:pullbacks_and_necessary}. One can find functional verb phrases connecting \fakebox{a man} to \fakebox{a woman}, for example ``has as mother". However, it is not straightforward to find a functional verb phrase that connects \fakebox{an object} to \fakebox{a number between 20 and 120} making the necessary diagram (analogous to (\ref{eqn:woman_number})) an endorsable equivalence between sentences. The reason, roughly, is that once one has regarded a man as an object, there is no aspect that applies to an arbitrary object which will correspond with the mother-having aspect of an arbitrary man. \lozend
\end{example}

The upshot of the above example is that linguistic structures give the necessary semantics to constrain mappings between ologs. This has important applications to databases, as we will show in Example~\ref{ex:constraints_help_databases}.


\subsection{Instantiated linguistic structures}\label{sec:instantiations}

In this section we study instantiated ologs and mappings between them. These are ologs for which each type has been assigned a set of examples. The authors of an olog should, and generally do, know more than just some types and relationships; they should also have in mind some examples of these types and relationships. For example, someone who writes an olog about dogs, say including the arrow $\fakebox{a dog}\To{\tn{has}}\fakebox{a name}$, probably knows some examples of dogs and their names. This information can be stored in an instantiation of the olog, which we will define in Definition~\ref{def:Instantiation}.

Instantiating an olog---filling it with conforming data---serves three purposes:
\begin{enumerate}[topsep=-1ex,itemsep=-1.5ex]
\item[\textbf{i.}] It gives users a place to store data about---examples of---their subject of interest.
\item[\textbf{ii.}] It validates the olog as a mathematical structure.
\item[\textbf{iii.}] It differentiates between different author sets who endorse the same conceptual scheme (we mean a noun phrase, a sentence, an equivalence between sentences, or an olog).
\end{enumerate}
The first of these purposes is probably the most important, but it is also straightforward, so we briefly explain the other two. Issues of functional connectivity and endorsed facts (see \Defcolor~\ref{def:FunctionRel}~and~\ref{def:DeclaredEq}) rely on the authors' understanding of mathematical functions and their compositions. By instantiating an olog, the users validate that understanding. 

Another reason to instantiate an olog is to differentiate one group of authors from another, even if they use the same conceptual scheme. For example, consider the following linguistic structure:
\begin{equation} \label{fig:fatherhood}
\hspace{-2.5cm} \mbox{$L$} = \parbox{1.3in}{\fbox{
\begin{tikzpicture}
\matrix (m) [ampersand replacement=\&, matrix of math nodes, row sep=0em, column sep=3em, text height=1.5ex, text depth=1.25ex]
{ \obox{a person} \& \obox{a father} \\ };
\path[->]
(m-1-1) edge node[above] {$\footnotesize\textrm{has}$} (m-1-2);
\end{tikzpicture} }}
\end{equation}

One author set may be interested in the fathers of US politicians (\eg George W. Bush's father is George H. W. Bush), whereas another set may be interested in the fathers of famous mathematicians (\eg Emmy Noether's father is Max Noether). The same linguistic expressions in an olog can house multiple instantiations. 

The mathematical motivation behind instantiations comes from the concept of set-valued functors $\cat{C} \longrightarrow \Set$ as database instances, introduced by \cite{Rosebrugh} and rediscovered by the second author in \cite{SpivakData}. The same notion was defined for ologs in \cite{SpivakOlogs}, where it was assumed that every olog comes equipped with such a functor $\cat{C}\too\Set$. We find three problems with this:
\begin{enumerate}[topsep=-1ex,itemsep=-1.5ex]
\item[\textbf{i.}] An olog can exist before one has recorded the corresponding examples.
\item[\textbf{ii.}] The examples should have something to do with the linguistics expressions set in the olog.
\item[\textbf{iii.}] Two authors may endorse the same linguistic expressions in an olog but have different examples.
\end{enumerate}
We have commented on \textbf{iii} above, and \textbf{i} is straightforward. We explain \textbf{ii}, which is probably the most important, in Remark~\ref{rem:instances_conform}.

\begin{remark}\label{rem:instances_conform}
In \cite{SpivakOlogs}, there was no assurance that the functor $I\colon\cat{C}\too\Set$ had anything to do with the linguistic expressions on $\cat{C}$. So a type $c=\fakebox{a dog}$ would be mapped to a set, but there was nothing ensuring that it was a set of dogs. Of course, such a thing cannot be ensured mathematically, but in Definition~\ref{def:validation} we did the next best thing, by providing a sentence for authors to endorse. \lozend
\end{remark}

\begin{example}
Consider the person-father olog $(\cat{C}, L)$ from \Figcolor~(\ref{fig:fatherhood}), reproduced here:
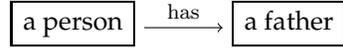
\[
\begin{tikzpicture}
\matrix (m) [ampersand replacement=\&, matrix of math nodes, row sep=0em, column sep=3em, text height=1.5ex, text depth=1.25ex]
{ \obox{a person} \& \obox{a father} \\ };
\path[->]
(m-1-1) edge node[above] {$\footnotesize\textrm{has}$} (m-1-2);
\end{tikzpicture} 
\]
An author $s$ might endorse that George W. Bush is a person, \ie $s\models (\mbox{George W. Bush}$ $\mbox{: a person})$, and similarly that George H. W. Bush is a father. Suppose, following (\ref{eqn:correspondence}) that $s$ also agrees with the sentence:
\begin{center}
``George W. Bush is a person, which has a father, namely \ul{George H. W. Bush.}"
\end{center}
Then $s$ endorses that ``has" is a correspondence between George W. Bush as a person and George H. W. Bush as a father. \lozend 
\end{example}

Let $U \colon \iEng \too \Eng$ and $\T \colon \iEng \too \Set$ be the functors from Proposition~\ref{prop:Instantiation}, and recall from Proposition~\ref{prop:Instantiation} that the induced functor $\iEng\inj\Eng\times\Set$ is an inclusion. The following definition is motivated to adapt the notion of instance functor given in \cite{SpivakOlogs} to our concept of linguistic structures.

\begin{definition}\label{def:Instantiation}
Let $\cat{C}$ be a category. An \textbf{instantiated linguistic structure} on $\cat{C}$ is a functor $\omfL\colon\cat{C}\to\iEng$. An \textbf{instantiated olog} is a pair $(\cat{C},\omfL)$, where $\omfL$ is an instantiated linguistic structure on $\cat{C}$.

Given a linguistic structure $L \colon \cat{C} \too \Eng$ and a set-valued functor $I \colon \cat{C} \too \Set$, we say that $I$ \textbf{conforms} to $L$ if $(L, I)$ factors through $\iEng$,
\begin{figure}[H]
\centering
\belowcaptionskip = -15pt
\begin{tikzpicture}[description/.style={fill=white,inner sep=2pt}]
\matrix (m) [matrix of math nodes, row sep=2.5em, column sep=1.75em, text height=1.75ex, text depth=0.25ex]
{ \cat{C} & & \Eng \times \Set \\ & \iEng \\ };
\path[->]
(m-1-1) edge node[above] {\footnotesize$(L, I)$} (m-1-3) edge node[below,sloped] {\footnotesize$\omfL$} (m-2-2);
\path[right hook->]
(m-2-2) edge node[below,sloped] {\footnotesize$(U, \T)$} (m-1-3);
\end{tikzpicture}
\caption{Functors $\cat{C} \too \Set$ conforming linguistic structures.}\label{fig:conforms}
\end{figure}
The lax functor $\omfL \colon \cat{C} \too \iEng$ will be sometimes referred as an \textbf{instantiation} of $L$. We may refer to $L$ and $I$ (respectively) as the \textbf{underlying linguistic structure} and \textbf{underlying instance} of $\omfL$. ~%
\lozend
\end{definition}

It is this notion, that an instance should conform to the linguistic structure, which we find missing in \cite{SpivakOlogs}; see Remark~\ref{rem:instances_conform}.

\begin{remark}\label{rem:instantiated_ling_char} \
If $\omfL \colon \cat{C} \too \iEng$ is a lax functor, then we obtain an instantiated linguistic structure $(L,I)$ on $\cat{C}$ by setting $L := U \circ \omfL$ and $I := \T \circ \omfL$; by definition the instance $I$ conforms to the linguistic structure $L$. We therefore do not distinguish between functors $\cat{C}\too\iEng$ and a pair $(L,I)$ of functors $L \colon \cat{C} \too \Eng$ and $I \colon \cat{C} \too \Set$, with $L$ lax, such that $I$ conforms to $L$. ~%
\lozend
\end{remark}

\begin{example}\label{ExampleBush} 
Recall again the person-father olog $(\cat{C}, L)$ from \Figcolor~\ref{fig:fatherhood}, reproduced here:
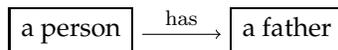
\[\begin{tikzpicture}
\matrix (m) [ampersand replacement=\&, matrix of math nodes, row sep=0em, column sep=3em, text height=1.5ex, text depth=1.25ex]
{ \obox{a person} \& \obox{a father} \\ };
\path[->]
(m-1-1) edge node[above] {$\footnotesize\textrm{has}$} (m-1-2);
\end{tikzpicture} 
\] 
We adopt a tabular description similar to that used in \cite{SpivakData}, except with column headings taken from (\ref{eqn:correspondence}). Using it, we can record the data of a functor $I\colon \cat{C} \longrightarrow \Set$ as follows: %
\footnote{For a sentence $(\N_1, \V, \N_2)$, the instance will fit into a table with $\qt{\N_1}$ as the head of the first column and $\qt{\V \N_2}$ as the head of the second column.}
\[
\begin{tabular}{ | l | l | }
\bhline
\textbf{a person} & \textbf{has a father, namely} \\ 
\bbhline
George W. Bush & George H. W. Bush \\
\hline 
Jeb Bush & George H. W. Bush\\
\hline
Emmy Noether & Max Noether \\
\bhline
\end{tabular} 
\hspace{.6in}
\begin{tabular}{ | l |  }
\bhline
\textbf{a father} \\ 
\bbhline
George H. W. Bush \\
\hline 
Max Noether\\
\hline
Bill Clinton\\
\bhline
\end{tabular} \]
This table then shows three correspondences, associated to the arrow labeled ``has":
\begin{itemize}[topsep=-1ex,itemsep=-1.5ex]
\item George W. Bush is a person, which has a father, namely \ul{George H. W. Bush};

\item Jeb Bush is a person, which has a father, namely \ul{George H. W. Bush}; and

\item Emmy Noether is a person, which has a father, namely \ul{Max Noether}. 
\end{itemize}
An author who endorses the six tokens and three correspondences shown here then also endorses $I$ as an instantiation of the olog $(\cat{C}, L)$. ~%
\lozend
\end{example}

\begin{remark}
The rules of instantiated English, as presented in Definition~\ref{def:type_aspect_fact_instance}, can be understood in terms of how people should interact with a database, as in \cite{SpivakData}. The database corresponding to an olog would have tables named by noun phrases and columns labeled by verb phrases. A person should only enter a new row in table $N$ if they understand (endorse) that noun phrase as a set, and thus can evaluate whether something was a member of it or not. A person should only fill a cell in column $V$ if they understand how it is a function, and thus can evaluate how the row corresponds to a row in some foreign table. {\color{white}{srsrd}}~%
\lozend
\end{remark}


\subsection{Instantiated functors}\label{sec:instantiated_functors}

We finish this section providing a notion of mappings between two instantiated ologs $(\cat{C}, \overline{L})$ and $(\cat{D}, \overline{M})$ in terms of lax transformations, as in Definition~\ref{def:LinguisticFunctor}, by constraining the set of linguistic functors $(\cat{C}, L) \longrightarrow (\cat{D}, M)$ to those that respect the tokens and their correspondences (See Definition~\ref{def:validation}). We will see in Remark~\ref{rem:natural_conforms} how this concept adapts the notion of \textit{meaningful functor} introduced in \cite[\S 4]{SpivakOlogsPre} to the context of linguistic functors.

\begin{definition}\label{def:InstantiatedFunctor}
Let $(\cat{C}, \omfL)$ and $(\cat{D}, \omfM)$ be instantiated ologs, as in Definition~\ref{def:Instantiation}. An \textbf{linguistic functor} between them, denoted $(F, F^{\natural}) \colon (\cat{C}, \omfL) \too (\cat{D}, \omfM)$ consists of a functor $F \colon \cat{C} \too \cat{D}$ together with a lax transformation $F^{\natural} \colon \omfL \leadsto \omfM \circ F$ as shown below:
\begin{figure}[H]
\belowcaptionskip = -15pt
\centering
\begin{tikzpicture}[description/.style={fill=white,inner sep=2pt}]
\matrix (m) [matrix of math nodes, row sep=0.25em, column sep=2.5em, text height=1.75ex, text depth=1.25ex]
{ \cat{C} & & \cat{D} \\ & \overset{F^{\natural}}{\leadsto} \\ & \iEng \\ };
\path[->]
(m-1-1) edge node[above] {\footnotesize$F$} (m-1-3) edge node[left] {\footnotesize$\omfL$\mbox{ \ }} (m-3-2)
(m-1-3) edge node[right] {\footnotesize\mbox{ \ }$\omfM$} (m-3-2);
\end{tikzpicture}
\caption{Instantiated functor.}\label{fig:instantiated_functor}
\end{figure}
If $F^{\natural}$ is the identity, in which case $\omfL = \omfM \circ F$, we say that the linguistic structure $\omfL$ is the \textbf{pullback of $\omfM$ along $F$} and write $\omfL = F^*(\omfM)$.

We define the \textbf{category of instantiated ologs}, denoted $\InstOlog$, to be the category whose objects are instantiated ologs and whose morphisms are instantiated linguistic functors. \lozend
\end{definition}

\begin{definition}\label{def:natural_conforms} 
Suppose that $\omfL, \omfM \colon \cat{C} \too \iEng$ are instantiated linguistic structures on $\cat{C}$, where $\omfL=(L,I)$ and $\omfM=(M,J)$ are the underlying linguistic structures and conforming instances. If $\alpha\colon L\leadsto M$ is a morphism of linguistic structures, \ie a lax transformation, and $p\colon I\Rightarrow J$ is a morphism of instances, \ie a natural transformation,
\begin{figure}[H]
\belowcaptionskip = -15pt
\centering
\begin{tikzpicture}[description/.style={fill=white,inner sep=2pt}]
\matrix (m) [matrix of math nodes, row sep=2.5em, column sep=10em, text height=1.75ex, text depth=1.25ex]
{ & \Eng \\ \cat{C} \\ & \Set \\ };
\path[->]
(m-2-1) edge [bend left=16] node[above] {\footnotesize$L$} (m-1-2) edge [bend right=20] node[below, near end] {\footnotesize$M$} node[above] {$\alpha \Downarrow\mbox{ \ \ \ }$} (m-1-2)
(m-2-1) edge [bend left=16] node[above, near end] {\footnotesize$I$} (m-3-2) edge [bend right=20] node[below] {\footnotesize$J$} node[above] {$p \Downarrow$} (m-3-2);
\end{tikzpicture}
\caption{Morphisms $\alpha$ of linguistic structures, and morphisms $p$ of instances.}\label{fig:mor_ling_inst}
\end{figure}
then we say that $p$ \textbf{conforms to} $\alpha$ if the lax transformation $(\alpha,p)\colon\cat{C}\to\Eng\times\Set$ factors through $\iEng\inj\Eng\times\Set$. 

Explicitly, for each object $c\in\cat{C}$, we have a verb phrase $\alpha_c\colon L(c) \leadsto M(c)$ and a function $p_c\colon I(c)\to J(c)$. We say that $p$ conforms to $\alpha$ if for each $c\in\cat{C}$ and token $x\in I(c)$, each author $s\in\Auth(\alpha_c)$ endorses the correspondence $\alpha_c(x,p_c(x))$. ~%
\lozend
\end{definition}

\begin{remark}\label{rem:natural_conforms} 
In \cite{SpivakOlogsPre}, a mapping between ologs is called a \textit{meaningful functor}; it is given by a functor $F \colon \cat{C} \too \cat{D}$ between two ologs $\cat{C}$ and $\cat{D}$, and a natural transformation $p\colon I \Rightarrow J \circ F$, where $I \colon \cat{C} \too \Set$ and $J \colon \cat{D} \too \Set$ are instances of $\cat{C}$ and $\cat{D}$. One particular issue with this notion is that there is no guarantee that the data defined by $p$ conforms to the linguistics defined by $F$ (nor that $I$ and $J$ conform to the linguistic structures on the olog). 

Our definition of instantiated functor remedies these issues. It is not hard to note from Definitions~\ref{def:InstantiatedFunctor} and \ref{def:natural_conforms}, along with Proposition~\ref{prop:Instantiation}, that every instantiated functor $(\cat{C}, \omfL) \too (\cat{D}, \omfM)$ is given by a functor $F \colon \cat{C} \too \cat{D}$, a lax transformation $F^\sharp \colon L \leadsto M \circ F$, and a natural transformation $F^\flat \colon I \Rightarrow J \circ F$ conforming to $F^\sharp$, where $L = U \circ \omfL$, $M = U \circ \omfM$, $I = \T \circ \omfL$ and $J = \T \circ \omfM$ as in Remark~\ref{rem:instantiated_ling_char}. This can be understood better with the help of Example~\ref{ex:constraints_help_databases}. ~%
\lozend
\end{remark}

In some situations, it is more natural to consider morphisms between instantiated ologs by reversing the lax transformations and natural transformations involved. Doing this yields another category of instantiated ologs, in which all the results obtained so far also hold. We explain this better in the following remark.

\begin{remark}
Suppose we are given an inclusion functor of categories $i \colon \cat{C} \too \cat{D}$, \ie $\cat{C}$ is a subcategory of $\cat{D}$, and that both are underlying categories of instantiated linguistic structures. Then there should be more people who understand $\cat{C}$, and therefore more data on $\cat{C}$ (than $\cat{D}$). For example, consider the following inclusion of ologs:
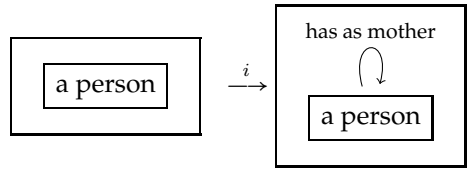
\begin{equation*}
\parbox{1.1in}{\fbox{ \begin{tikzpicture}[description/.style={fill=white,inner sep=2pt}] 
\matrix (m) [ampersand replacement=\&, matrix of math nodes, row sep=3.5em, column sep=0em] 
{ \obox{a person} \\ }; 
\path[->] 
; 
\end{tikzpicture} } 
}
\xrightarrow{\;\;i\;\;}
\parbox{1.5in}{\fbox{ \begin{tikzpicture}[description/.style={fill=white,inner sep=2pt}] 
\matrix (m) [ampersand replacement=\&, matrix of math nodes, row sep=3.5em, column sep=0em] 
{ \obox{a person} \\ }; 
\path[->] 
(m-1-1) edge [loop above] node {\footnotesize$\mbox{has as mother}$} (m-2-1); 
\end{tikzpicture} }
}
\end{equation*}  
So in this case, there would be a map of instances $q\colon J \circ i \Rightarrow I$, where $I$ is an instance on $\cat{C}$ and $J$ is an instance on $\cat{D}$. There would also be a map of linguistic structures $\beta\colon M \circ i \leadsto L$, where $L$ is a linguistic structure on $\cat{C}$ and $M$ is a linguistic structure on $\cat{D}$, to which $q$ conforms. In fact, note that if $J \circ i \Rightarrow I$ is the direction of instances, then it is only natural that the linguistics could go ``the other way" $M \circ i \leadsto L$:

\begin{equation}\label{eqn:co_instantiated}
\parbox{1.3in}{
\begin{tikzpicture}[description/.style={fill=white,inner sep=2pt}]
\matrix (m) [ampersand replacement=\&, matrix of math nodes, row sep=0.4em, column sep=2em, text height=1.75ex, text depth=1.25ex]
{ \cat{C} \& \& \cat{D} \\ \& \overset{(\beta,q)}{\Longleftarrow} \\ \& \iEng \\ };
\path[->]
(m-1-1) edge node[above] {\footnotesize$i$} (m-1-3) edge node[left] {\footnotesize$(L,I)$} (m-3-2)
(m-1-3) edge node[right] {\footnotesize$(M,J)$} (m-3-2);
\end{tikzpicture}
}
\end{equation}
We define the category $\InstOlog^\dagger$ of co-instantiated ologs to have the same objects as $\InstOlog$, but where morphisms are given as in \eqref{eqn:co_instantiated}.\footnote{The $\dagger$ symbol in our case has nothing to do with the dagger category \cite{WikiDagger}.} ~%
\lozend
\end{remark}

\begin{remark}\label{rem:Grothendieck2} 
As we did with the category $\Olog$ in Remark~\ref{rem:Grothendieck}, we note there is an evident functor $\InstOlog \to \Cat$, sending every instantiated olog $(\cat{C}, \omfL)$ to its underlying category $\cat{C}$. This is a Grothendieck fibration whose cartesian arrows coincide with the \textit{strongly meaningful functors} defined in \cite[\S 4]{SpivakOlogsPre}.  ~%
\lozend
\end{remark}

\begin{example}\label{ex:constraints_help_databases}
Suppose we are given two categories $\cat{C},\cat{D}$, which we think of as database schemas (as in \cite{SpivakData}), a functor $F \colon \cat{C} \too \cat{D}$, and two instantiations $I\colon\cat{C}\too\Set$ and $J\colon\cat{D}\too\Set$. Suppose that these two databases are to be merged. We are asked to find a morphism $F^\flat \colon I \Rightarrow J \circ F$. In this example, we show that our job will be easier if $F$ has been equipped with a morphism $F^{\#}$ of linguistic structures, so we can find $F^\flat$ conforming to $F^\#$.

In order to emphasize the issue, we suppose that $\cat{C}$ and $\cat{D}$ are both single-object categories 
$$\cat{C}=\fbox{$\bullet^1$} \mbox{ \ and \ } \cat{D}=\fbox{$\bullet^a$}$$ 
equipped with the following linguistic structures: 
\begin{align*}
\mbox{ \ \ } (\cat{C}, L) & = \parbox{1.3in}{\fbox{
\begin{tikzpicture}[description/.style={fill=white,inner sep=2pt}] 
\matrix (m) [ampersand replacement=\&, matrix of math nodes, row sep=5em]
{ \ovbox{1}{a human} \\ }; 
\end{tikzpicture} }} \mbox{ and \ \ \ \ \ \ } 
(\cat{D}, M) = \parbox{1.3in}{\fbox{ 
\begin{tikzpicture}[description/.style={fill=white,inner sep=2pt}] 
\matrix (m) [ampersand replacement=\&, matrix of math nodes, row sep=5em]
{ \ovbox{a}{a person} \\ }; 
\end{tikzpicture} }} 
\end{align*} 
Suppose we want to compare two instantiations, $I\colon\cat{C}\too\Set$ and $J\colon\cat{D}\too\Set$, which are represented by the following tables:
\[
I(1):=\begin{tabular}{ | c |}
\bhline
\textbf{a human} \\ 
\bbhline
Emmy Noether \\
\hline 
George W. Bush \\
\bhline
\end{tabular}
\hspace{.8in}
J(a):=\begin{tabular}{ | c |}
\bhline
\textbf{a person} \\ 
\bbhline
Emmy Noether \\
\hline 
Max Noether \\
\hline
Bill Clinton\\
\hline
George H. W. Bush\\
\hline
George W. Bush \\
\bhline
\end{tabular}
\]

To compare instances on different schemas, we first need a functor between them. In our case there is a unique functor $F\colon\cat{C} \too \cat{D}$ (it sends $1\mapsto a$), so this is not an issue. With this functor in hand, we can pull back $J$ to an instantiation $J \circ F$ on $\cat{C}$, and attempt to compare it to $I$. 

The purpose of this example is to show that the choice of linguistic structure $F^{\#}$ on $F$ is an important aid to making this comparison, \ie to choosing a database homomorphism $I\Rightarrow J \circ F$ out of the $5^2=25$ possible choices. Consider the following two morphisms for $F$:
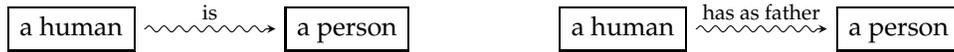
\[ \begin{tikzpicture}[description/.style={fill=white,inner sep=2pt}] 
\matrix (m) [matrix of math nodes, row sep=3em, column sep=5em]
{ \obox{a human} & \obox{a person} & \obox{a human} & \obox{a person} \\ };
\path[->]
(m-1-1) edge [snake it] node[above] {\footnotesize\mbox{is}} (m-1-2)
(m-1-3) edge [snake it] node[above] {\footnotesize\mbox{has as father}} (m-1-4);
\end{tikzpicture} \] 
We denote the first by $\alpha$ and the second by $\beta$, \ie $\alpha,\beta \colon L \leadsto M \circ F$. Clearly, these morphisms give a useful hint at the intended semantics for the mapping. The only natural transformation $I\Rightarrow J \circ F$ we endorse as conforming to $\alpha$ (see Definition~\ref{def:natural_conforms}) is $p$, and the only one that conforms to $\beta$ is $q$, as shown below:
\[ (\alpha,p) :=\;\; 
\begin{tabular}{ | c | c | }
\bhline
\textbf{a human} & \textbf{is a person, namely} \\ 
\bbhline
Emmy Noether & Emmy Noether \\
\hline 
George W. Bush & George W. Bush \\
\bhline
\end{tabular} \]
\[ (\beta,q) := \;\;
\begin{tabular}{ | c | c | }
\bhline
\textbf{a human} & \textbf{has as father a person, namely} \\ 
\bbhline
Emmy Noether & Max Noether \\
\hline 
George W. Bush & George H. W. Bush \\
\bhline
\end{tabular} \]

This example complements, but also extends the scope of ``meaningful functors" found in \cite[\S 4]{SpivakOlogsPre}. \lozend
\end{example}


\section*{Acknowledgements}

The authors are supported by the following grants: Office of Naval Research ONR N00014131 0260, Air Force Office of Scientific Research AFOSR FA9550-14-1-0031, and National Aeronautics and Space Administration NASA (Langley Research Center) NNH13ZEA001N-SSAT. 

The authors want to thank Dr. Patrick Schultz from MIT for patiently listening to some of the authors' discussions on this research, and for his comments which cleared up some doubts at the abstract level of this paper. 

Finally, special thanks to Dr. Spencer Breiner from NIST, who patiently read the first manuscript, and whose comments and corrections improved the quality of the present paper. The first author also thanks him for the discussions both had during his visit to MIT in April 2015, which inspired alternative ways to present some of the most important concepts in the theory of ologs.


\end{document}